\documentclass[12pt]{amsart}
\usepackage{amscd,amsmath,amssymb,amsfonts}
\usepackage[all]{xypic}
\theoremstyle{plain}
\newtheorem{thm}[subsection]{Theorem}
\newtheorem{lem}[subsection]{Lemma}
\newtheorem{cor}[subsection]{Corollary}
\newtheorem{prop}[subsection]{Proposition}

\theoremstyle{definition}
\newtheorem{defn}[subsection]{Definition}

\newtheorem{rmks}[subsection]{Remarks}

\numberwithin{subsection}{section}
\numberwithin{equation}{section}

\newcommand{\eq}[2]{\begin{equation}\label{#1}#2 \end{equation}}
\newcommand{\ml}[2]{\begin{multline}\label{#1}#2 \end{multline}}
\newcommand{\ga}[2]{\begin{gather}\label{#1}#2 \end{gather}}
\newcommand{\gr}{{\rm gr}}

\newcommand{\surj}{\twoheadrightarrow}
\newcommand{\inj}{\hookrightarrow}
\newcommand{\red}{{\rm red}}

\newcommand{\Hom}{{\rm Hom}}

\newcommand{\Spec}{{\rm Spec \,}}

\newcommand{\Char}{{\rm char}}

\newcommand{\Gal}{{\rm Gal}}

\newcommand{\sA}{{\mathcal A}}

\newcommand{\sC}{{\mathcal C}}
\newcommand{\sD}{{\mathcal D}}

\newcommand{\sH}{{\mathcal H}}
\newcommand{\sI}{{\mathcal I}}
\newcommand{\sJ}{{\mathcal J}}
\newcommand{\sK}{{\mathcal K}}
\newcommand{\sL}{{\mathcal L}}

\newcommand{\sN}{{\mathcal N}}
\newcommand{\sO}{{\mathcal O}}
\newcommand{\sP}{{\mathcal P}}

\newcommand{\sR}{{\mathcal R}}

\newcommand{\C}{{\mathbb C}}

\newcommand{\F}{{\mathbb F}}

\renewcommand{\H}{{\mathbb H}}

\newcommand{\N}{{\mathbb N}}
\renewcommand{\P}{{\mathbb P}}
\newcommand{\Q}{{\mathbb Q}}
\newcommand{\R}{{\mathbb R}}

\newcommand{\U}{{\mathbb U}}

\newcommand{\W}{{\mathbb W}}

\newcommand{\Z}{{\mathbb Z}}

\newcounter{alphab}[subsection]
\newcommand{\alphab}
{\stepcounter{alphab}
\noindent\makebox[0.5cm]{}{\normalfont\alph{alphab})}\hspace{0.2cm}}

\newcounter{romain}[subsection]
\newcommand{\romain}
{\stepcounter{romain}
\noindent\makebox[1.5cm][r]{{\normalfont(\roman{romain})}\hspace{0.3cm}}}

\newcommand{\cf}{\textit{cf.\ }}
\newcommand{\ie}{\textit{i.e.\ }}

\newcommand{\Frac}{\mathrm{Frac}}
\newcommand{\Ker}{\mathrm{Ker}}
\newcommand{\Coker}{\mathrm{Coker}}
\renewcommand{\Im}{\mathrm{Im}}

\newcommand{\Spf}{\mathrm{Spf}}
\newcommand{\Spm}{\mathrm{Spm}}
\newcommand{\spm}{\mathrm{sp}}
\newcommand{\cpt}{_{\mathrm{c}}}

\newcommand{\cris}{_{\mathrm{crys}}}

\newcommand{\rig}{_{\mathrm{rig}}}
\newcommand{\rigc}{_{\mathrm{rig,c}}}
\newcommand{\tot}{_{\mathrm{t}}}
\newcommand{\Omd}[1]{\Omega^{^{_{_{\bullet}}}}_{#1}}

\newcommand{\WOmd}[2]{W_{#1}\Omega^{^{_{_{\bullet}}}}_{#2}}
\newcommand{\Db}{D^{\mathrm{b}}}
\newcommand{\sPh}{\widehat{\sP}}

\newcommand{\bul}{^{_{_{\bullet}}}\!}
\newcommand{\hbul}{^{^{_{_{\bullet}}}}}

\newcommand{\lbul}{_{_{^\bullet}}}

\newcommand{\ua}{\underline{a}}
\newcommand{\sk}{\mathrm{sk}}
\newcommand{\cosk}{\mathrm{cosk}}
\newcommand{\kbar}{\bar{k}}
\newcommand{\Ebar}{\,\overline{\!E}}
\newcommand{\Fbar}{\,\overline{\!F}}
\newcommand{\Kbar}{\,\overline{\!K}}
\newcommand{\Ubar}{\overline{U}}
\newcommand{\Ob}{\mathrm{Ob}}

\begin{document}

\title[Witt Vector Cohomology of Singular Varieties]{
On Witt Vector Cohomology for Singular Varieties}
\author{Pierre Berthelot}
\address{IRMAR, Universit\'e de Rennes 1,
Campus de Beaulieu,
35042 Rennes cedex, France
}
\email{pierre.berthelot@univ-rennes1.fr}
\author{Spencer Bloch}
\address{Dept. of Mathematics,
University of Chicago,
Chicago, IL 60637,
USA}
\email{bloch@math.uchicago.edu}
\author{H\'el\`ene Esnault}
\address{Mathematik,
Universit\"at Duisburg-Essen, FB6, Mathematik, 45117 Essen, Germany}
\email{esnault@uni-essen.de}
\thanks{The first author has been supported by the research network
{\it Arithmetic Algebraic Geometry} of the European Community
(Contract MRTN-CT-2003-504917).}

\begin{abstract}
Over a perfect field $k$ of characteristic $p > 0$, we construct a
``Witt vector cohomology with compact supports'' for separated
$k$-schemes of finite type, extending (after tensorisation with $\Q$)
the classical theory for proper $k$-schemes. We define a canonical
morphism from rigid cohomology with compact supports to Witt vector
cohomology with compact supports, and we prove that it provides an
identification between the latter and the slope $< 1$ part of the
former. Over a finite field, this allows one to compute congruences for
the number of rational points in special examples. In particular, the
congruence modulo the cardinality of the finite field of the number of
rational points of a theta divisor on an abelian variety does not
depend on the choice of the theta divisor. This answers positively a
question by J.-P. Serre.
\end{abstract}

\maketitle
\thispagestyle{empty}

\setcounter{tocdepth}{1}
\tableofcontents

\section{Introduction}
\label{Introduction}
\medskip

Let $k$ be a perfect field of characteristic $p > 0$, $W = W(k)$, $K =
\Frac(W)$. If $X$ is a proper and smooth variety defined over $k$, the
theory of the de Rham-Witt complex and the degeneration of the slope
spectral sequence provide a functorial isomorphism (\cite[III,
3.5]{B}, \cite[II, 3.5]{I})
\eq{sl1} {H^{\ast}\cris(X/K)^{<1}\ \xrightarrow{\ \sim\ \,}\ H^{\ast}(X, 
W\sO_X)_K,}
where $H^{\ast}\cris(X/K)^{<1}$ is the maximal subspace on which
Frobenius acts with slopes $< 1$, and the subscript $_K$ denotes 
tensorisation with $K$.

If we only assume that $X$ is proper, but maybe singular, we can use
rigid cohomology to generalize crystalline cohomology while retaining
all the standard properties of a topological cohomology theory. Thus,
the left hand side of $\eqref{sl1}$ remains well defined, and, when
$k$ is finite, the alternated product of the corresponding
characteristic polynomials of Frobenius can be interpreted as the
factor of the zeta function $\zeta(X,t)$ ``of slopes $< 1$'' (see
\ref{Pless1} for a precise definition). On the other hand, the
classical theory of the de Rham-Witt complex can no longer be directly
applied to $X$, but the sheaf of Witt vectors $W\sO_X$ is still
available. Thus the right hand side of $\eqref{sl1}$ remains also well
defined. As in the smooth case, this is a finitely generated
$K$-vector space endowed with a Frobenius action with slopes in
$[0,1[$ (Proposition \ref{finite}), which has the advantage of being
directly related to the coherent cohomology of $X$. It is therefore of
interest to know whether, when $X$ is singular, \eqref{sl1} can be 
generalized as an isomorphism
\ga{sl2} {H^{\ast}\rig(X/K)^{<1}\ \xrightarrow{\ \sim\ \,}\ H^{\ast}(X, 
W\sO_X)_K,}
where $H^{\ast}\rig(X/K)^{<1}$ denotes the subspace of slope $< 1$ 
of $H^{\ast}\rig(X/K)$.

Our main result gives a positive answer to this question. More
generally, we show that a ``Witt vector cohomology with compact
supports'' can be defined for separated $k$-schemes of finite type,
giving cohomology spaces $H^{\ast}\cpt(X, W\sO_{X,K})$ which are
finite dimensional $K$-vector spaces, endowed with a Frobenius action
with slopes in $[0,1[$. Then, for any such scheme $X$, the slope $< 1$
subspace of the rigid cohomology with compact supports of $X$ has the
following description:

\begin{thm}\label{SlopeHc}
Let $k$ be a perfect field of characteristic $p > 0$, $X$ a separated
$k$-scheme of finite type. There exists a functorial isomorphism
\ga{sl3} {H^{\ast}\rigc(X/K)^{<1}\ \xrightarrow{\ \sim\ \,}\
H^{\ast}\cpt(X, W\sO_{X,K}).}
\end{thm}

This is a striking confirmation of Serre's intuition \cite{Se} about
the relation between topological and Witt vector cohomologies. On the
other hand, this result bears some analogy with Hodge theory. Recall
from \cite{BH} (or \cite{De1}) that if $X$ is a proper scheme defined
over $\C$, then its Betti cohomology $H^{\ast}(X, \C)$ is a direct
summand of its de Rham cohomology $\H^{\ast}(X, \Omega^\bullet_X)$.
Coherent cohomology $H^{\ast}(X, \sO_X)$ does not exactly compute the
corner piece of the Hodge filtration, as would be an exact analogy
with the formulation of Theorem \ref{SlopeHc}, but it gives an upper
bound. Indeed, by \cite{E}, Proposition 1.2, and \cite{E2}, Proof of
Theorem 1.1, one has a functorial surjective map $H^{\ast}(X,
\sO_X)\surj {\rm gr}_F^0 H^{\ast}(X, \C)$.
\medskip

The construction of these Witt vector cohomology spaces is given in
section \ref{WittCohCompSupp}. If $U$ is a separated $k$-scheme of
finite type, $U \inj X$ an open immersion in a proper $k$-scheme, and
$\sI \subset \sO_X$ any coherent ideal such that $V(\sI) = X \setminus
U$, we define $W\sI = \Ker(W\sO_X \to W(\sO_X/\sI))$, and we show that
the cohomology spaces $H^{\ast}(X, W\sI_K)$ actually depend only on
$U$ (Theorem \ref{IndComp}). This results from an extension to Witt
vectors of Deligne's results on the independence on the
compactification for the construction of the $f_!$ functor for
coherent sheaves \cite{Ha1}.

When $U$ varies, these spaces are contravariant functors with respect
to proper maps, and covariant functors with respect to open
immersions. In particular, they give rise to the usual long exact
sequence relating the cohomologies of $U$, of an open subset $V
\subset U$, and of the complement $T$ of $V$ in $U$. Thus, they can be
viewed as providing a notion of Witt vector cohomology spaces with
compact support for $U$.  We define
$$H^{\ast}\cpt(U, W\sO_{U,K}):=H^{\ast}(X, W\sI_K).$$  
Among other properties, we prove in section \ref{SectDesc} that these
cohomology spaces satisfy a particular case of cohomological descent
(Theorem \ref{descent}) which will be used in the proof of Theorem
\ref{SlopeHc}.

The construction of the canonical homomorphism between rigid and Witt
vector cohomologies is given in section \ref{RigandWitt} (Theorem
\ref{ConstRigtoWitt}). First, we recall how to compute rigid
cohomology for a proper $k$-scheme $X$ when there exists a closed
immersion of $X$ in a smooth formal $W$-scheme $\P$, using the de Rham
complex of $\P$ with coefficients in an appropriate sheaf of
$\sO_{\P}$-algebras $\sA_{X,\P}$. When $\P$ can be endowed with a
lifting of the absolute Frobenius endomorphism of its special fibre, a
simple construction (based on an idea of Illusie \cite{I}) provides a
functorial morphism from this de Rham complex to $W\sO_{X,K}$. Using
\v{C}ech resolutions, this morphism can still be defined in the
general case as a morphism in the derived category of sheaves of
$K$-vector spaces. Then we obtain \eqref{sl3} by taking the morphism
induced on cohomology and restricting to the slope $<1$ subspace. By
means of simplicial resolutions based on de Jong's theorem, we prove
in section \ref{ProofMainTh} that it is an isomorphism. The proof
proceeds by reduction to the case of proper and smooth schemes, using
Theorem \ref{descent} and the descent properties of rigid cohomology
proved by Chiarellotto and Tsuzuki (\cite{CT}, \cite{Ts1},
\cite{Ts2}).
\medskip

We now list some applications of Theorem \ref{SlopeHc}, which are
developed in section \ref{Applications}. We first
remark that it implies a vanishing theorem. If $X$ is smooth
projective, and $Y\subset X$ is a divisor so that $U=X\setminus Y$ is
affine, then Serre vanishing says that $H^i(X, \sO_X(-nY))=0$ for $n$
large and $i<\dim(X)$. If $k$ has characteristic $0$, then we can take
$n=1$ by Kodaira vanishing.  One then has $H^i(X, \sO_X(-Y))=0$ for
$i<\dim(X)$. It is tempting to view the following corollary as an
analogue when $k$ has characteristic $p$\,:

\begin{cor} \label{Vanish}
Let $X$ be a proper scheme defined over a perfect field $k$ of
characteristic $p > 0$. Let $\sI \subset \sO_X$ be a coherent ideal
defining a closed subscheme $Y\subset X$ such that $U=X\setminus Y$
is affine, smooth and equidimensional of dimension $n$. Then 
$$ H^i_c(U, W\sO_{U,K})=H^i(X,
W\sI)_K = H^i(X, W\sI^r)_K=0$$ 
for all $i \neq n$ and all $r \geq 1$.
\end{cor}

Indeed, when $U$ is smooth and affine, $H^i\rig(U/K)$ can be
identified with Monsky-Washnitzer cohomology \cite{Be1}, and therefore
vanishes for $i > n = \dim(U)$. If moreover $U$ is equidimensional, it
follows by Poincar\'e duality \cite{Be2} that $H^i\rigc(U/K)=0$ for $i
< n$. So \eqref{sl3} implies that $H^i(X, W\sI)_K = 0$ for $i < n$. On
the other hand, the closure $\Ubar$ of $U$ in $X$ has dimension $n$,
and $H^i(X, W\sI)_K = H^i\cpt(U, W\sO_{U,K})$ does not change up to
canonical isomorphism if we replace $X$ by $\Ubar$. Therefore $H^i(X,
W\sI)_K = 0$ for $i > n$.

\medskip
When $k$ is a finite field $\F_q$, with $q = p^a$, Theorem
\ref{SlopeHc} implies a statement about zeta functions. By
(\cite{ELS}, Th\'eor\`eme II), the Lefschetz trace formula provides an
expression of the $\zeta$-function of $X$ as the alternating product
$$\zeta(X,t)=\prod_i P_i(X, t)^{(-1)^{i+1}},$$
where $P_i(X, t) = \det(1-t\phi|H^i\rigc(X/K))$, and $\phi = F^a$
denotes the $\F_q$-linear Frobenius endomorphism of $X$. On the other
hand, we define
\ga{I6}{P_i^W(X,t) = \det(1-t\phi|H^i\cpt(X, W\sO_{X,K})),\\
 \zeta^W(X,t)=\prod_i P_i^W(X,t)^{(-1)^{i+1}}.\notag}
If we denote by $\zeta^{<1}(X,t)$ the ``slope $< 1$ factor'' of
$\zeta(X,t)$ (\cf \ref{Pless1}), we get formally from \eqref{sl3}:

\begin{cor}\label{zetaW}
Let $X$ be a separated scheme of finite type over a finite field. Then
one has
\ga{fZetaW}{ \zeta^{<1}(X,t) = \zeta^W(X,t). }
\end{cor}

This result can be used to prove congruences mod $q$ on the number of
$\F_q$-rational points of certain varieties. The following theorem
answers a question of Serre, and was the initial motivation for this
work. Recall that an effective divisor $D$ on an abelian variety $A$
is called a theta divisor if $\sO_A(D)$ is ample and defines a 
principal polarization.

\begin{thm} \label{ConjSerre}
Let $\Theta, \Theta'$ be two theta divisors
on an abelian variety defined over a finite field
$\F_q$. Then:
\eq{CongSerre}{ |\Theta(\F_q)| \equiv |\Theta'(\F_q)| \mod q. }
\end{thm}

Actually, Serre's original formulation predicts that, on an abelian
variety defined over a field, the difference of the motives of
$\Theta$ and $\Theta'$ is divisible by the Lefschetz motive. Our
Theorem \ref{ConjSerre} answers the point counting consequence of it.

We also have more elementary point counting consequences. 

\begin{cor}[Ax \cite{Ax}, Katz \cite{Ka1}]\label{ThAx}
Let $D_1,\ldots,D_r \subset \P^n$ be hypersurfaces of degrees
$d_1,\ldots,d_r$, defined over the finite field $\F_q$. Assume that
$\sum_j d_j \leq n$. Then
\eq{CongAx}{ |(D_1\cap \ldots \cap D_r)(\F_q)|\equiv 1 \mod q. }
\end{cor}

We observe here that this congruence is the best approximation of the
results of Ax and Katz that can be obtained using Witt vector
cohomology, since this method only provides information on the slope
$< 1$ factor of the zeta function. It would be worthwhile to have for
higher slopes results similar to Theorem \ref{SlopeHc} which might
give the full Ax-Katz congruences.

We also remark that Ax's theorem has a motivic proof \cite{BEL}, which
of course is more powerful than this slope proof. Yet it is of
interest to remark that Theorem \ref{SlopeHc} applies here as well. As
for Theorem \ref{ConjSerre}, it seems more difficult to formulate a
motivic proof, as it would have to deal with non-effective motives
(see the discussion in \ref{RemonSerre}).

Let us mention finally the following general consequence of Theorem
\ref{SlopeHc}, which for example can be applied in the context of the
work of Fu and Wan on mirror congruences for Calabi-Yau varieties
(\cite{Wa}, \cite{WF}):

\begin{cor}\label{IsomOX}
Let $f : X \to Y$ be a morphism between two proper $\F_q$-schemes. If 
$f$ induces isomorphisms $f^{\ast} : H^i(Y, \sO_Y) \xrightarrow{\ 
\sim\ \,} H^i(X, \sO_X)$ for all $i \geq 0$, then
\eq{CongIsomOX}{ |X(\F_q)| \equiv |Y(\F_q)| \mod q. }
\end{cor}

\bigskip
\noindent \textit{Acknowledgements:} It is a pleasure to thank
Jean-Pierre Serre for his strong encouragement and his help. Theorem
\ref{ConjSerre} was the main motivation for this work. We  thank Luc
Illusie for useful discussions. The third named author thanks Eckart Viehweg for his interest and encouragement. Corollary \ref{Vanish} is inspired by  the analogy  on the relation between de Rham and coherent cohomology, which has been developed jointly with him.

\bigskip
\noindent \textit{Notations and conventions:} Throughout this article,
$k$ denotes a perfect field of characteristic $p$, $W = W(k)$, $\sigma
: W \to W$ the Frobenius automorphism of $W$, $K = \Frac(W)$. The
subscript $_K$ will denote tensorisation with $K$ over $W$. We recall
that, on any noetherian topological space, taking cohomology commutes
with tensorisation by $\Q$. Therefore the subscript $_K$ will be moved
inside or outside cohomology or direct images without further
justification.

We denote by $\Db(K)$ (resp.\ $\Db(X,K)$) the derived category of
bound\-ed complexes of $K$-vector spaces (resp.\ complexes of sheaves 
of $K$-vector spaces on a topological space $X$).

All formal schemes considered in this article are $W$-formal schemes 
for the $p$-adic topology.

\bigskip
\section{Witt vector cohomology with compact supports} 
\label{WittCohCompSupp}
\medskip

We give in this section some properties of Witt vector cohomology
which are a strong indication of its topological nature, and will be
used later in our proof of Theorem \ref{SlopeHc}. In particular, we
show how to attach Witt vector cohomology spaces with compact supports
to any separated $k$-scheme of finite type.

If $X$ is a scheme, $\sA$ a sheaf of rings on $X$, and $n \geq 1$, we
denote by $W_n\sA$, or $W_n(\sA)$ if confusion may arise, the sheaf of
Witt vectors of length $n$ with coefficients in $\sA$, and by $W\sA
= \varprojlim_n W_n\sA$, or $W(\sA)$, the sheaf of Witt vectors of
infinite length. If $\sI \subset \sA$ is an ideal, we denote by
$W_n\sI = \Ker(W_n\sA \to W_n(\sA/\sI))$, or $W_n(\sI)$, the sheaf of
Witt vectors $(a_0,a_1, \ldots, a_{n-1})$ such that $a_i$ is a section
of $\sI$ for all $i$, and we define similarly $W\sI$. Note that, when
$\sI$ is quasi-coherent, the canonical morphism $W\sI \to
\R\varprojlim_n W_n\sI$ is an isomorphism, as $H^1(U,W_n\sI)=0$ for
any affine open subset $U \subset X$ and any $n$, and the projective
system $\Gamma(U, W_n\sI)$ has surjective transition maps.

For any $X$, any sheaf of rings $\sA$ on $X$, and any ideal $\sI
\subset \sA$, we use the notations $\R\Gamma(X, W\sI)$ and
$H^{\ast}(X, W\sI)$ to denote the Zariski cohomology of the sheaf
$W\sI$. Thus, when $\sI$ is quasi-coherent, the canonical morphism
$$\R\Gamma(X, W\sI) \to \R\varprojlim_n \R\Gamma(X,W_n\sI)$$
is an isomorphism. When $X$ is a proper $k$-scheme and $\sI \subset
\sO_X$ is a coherent ideal, the cohomology modules $H^{\ast}(X,
W_n\sI)$ are artinian $W$-modules; then it follows from the
Mittag-Leffler criterium that the morphism
$$H^{\ast}(X, W\sI) \to \varprojlim_n H^{\ast}(X, W_n\sI)$$
is an isomorphism.

We will shorten notations by writing $W\sO_{X,K}$, $W\sI_K$ for 
$(W\sO_X)_K$, $(W\sI)_K$. We recall again that, when $X$ is noetherian, 
there is a canonical isomorphism
$$ \R\Gamma(X, W\sI)_K \xrightarrow{\ \sim\ \,} \R\Gamma(X, W\sI_K). $$ 

In this article, we will be particularly interested in the $K$-vector
spaces $H^{\ast}(X, W\sO_{X,K})$, and in their generalizations
$H^{\ast}(X, W\sI_K)$. Our main observation is that, in contrast to
the $W$-modules $H^{\ast}(X, W\sO_X)$, which are sensitive to
nilpotent sections of $\sO_X$, they behave like a topological
cohomology theory for separated $k$-schemes of finite type. The
following easy proposition is somehow a key point.

\begin{prop}\label{TopInv} Let $X$ be a $k$-scheme of finite type.\\
\romain The canonical homomorphism
$$ W\sO_{X,K} \to W\sO_{X_{\red},K}$$
is an isomorphism, and induces a functorial isomorphism
$$ H^*(X,W\sO_{X,K}) \xrightarrow{\ \sim\ \,}
H^*(X_{\red},W\sO_{X_{\red},K}),$$
compatible with the action of $F$ and $V$.\\
\romain Let $\sI, \sJ \subset \sO_X$ be coherent ideals, and assume
$\sqrt{\sI} = \sqrt{\sJ}$, i.e. there exists $N\ge 1$ such that $\sI^N
\subset \sJ$ and $\sJ^N \subset \sI$. Then there is a canonical
identification 
$$ W\sI_K \cong W\sJ_K, $$
inducing a functorial isomorphism
$$ H^*(X, W\sI_K) \cong H^*(X,W\sJ_K), $$
compatible with the action of $F$ and $V$.
\end{prop}

\begin{proof} 
Let $\sN = \Ker(\sO_X \to \sO_{X_{\red}})$. To prove the first claim
of assertion (i), it suffices to show that $W\sN_K = 0$. But the
action of $p$ is invertible, and $p=VF=FV$, so it suffices to show $F$
is nilpotent. This is clear since $F$ acts on $W\sN$ by raising
coordinates to the $p$-th power. Taking cohomology, the second claim
follows.

Assertion (ii) follows from assertion (i) applied to $\sO_X/\sI$ and 
$\sO_X/\sJ$. 
\end{proof}

We now prove the existence of Mayer-Vietoris exact sequences for
Witt vector cohomology.

\begin{prop}\label{MayerViet}
Let $X$ be a $k$-scheme of finite type, $X_1, X_2 \subset X$ two 
closed subschemes such that $X = X_1 \cup X_2$, and $Z = X_1 \cap X_2$. 
There is an exact sequence
\eq{MVSeq}{ 0 \to W\sO_{X,K} \to W\sO_{X_1,K} \oplus W\sO_{X_2,K} \to
W\sO_{Z,K} \to 0, }
providing a Mayer-Vietoris long exact sequence
\ml{}{ \cdots \to H^i(X, W\sO_{X,K}) \to H^i(X_1, W\sO_{X_1,K}) \oplus 
H^i(X_2, W\sO_{X_2,K}) \notag\\
\to H^i(Z, W\sO_{Z,K}) \to H^{i+1}(X, W\sO_{X,K}) \to \cdots .\notag }
\end{prop}

\begin{proof}
Let $\sI_1, \sI_2$ be the ideals of $\sO_X$ defining $X_1$ and $X_2$.
Thanks to \ref{TopInv}, we may assume that $\sO_Z =
\sO_X/(\sI_1+\sI_2)$. It is easy to check that
$W(\sI_1+\sI_2)=W(\sI_1)+W(\sI_2)$. From the exact sequence
$$ 0 \to \sO_X/(\sI_1\cap\sI_2) \to \sO_X/\sI_1\oplus\sO_X/\sI_2 \to 
\sO_X/(\sI_1+\sI_2) \to 0,$$
we can then deduce an exact sequence
$$ 0 \to W(\sO_X/(\sI_1\cap\sI_2)) \to W\sO_{X_1}\oplus 
W\sO_{X_2} \to W\sO_Z \to 0.$$
Since $X = X_1 \cup X_2$, $\sI_1 \cap \sI_2$ is a nilpotent ideal, 
and Proposition \ref{TopInv} implies that $W\sO_{X,K} 
\xrightarrow{\ \sim\ \,} W(\sO_X/(\sI_1\cap\sI_2))_K$. This gives 
the short exact sequence of the statement. The long one follows by
taking cohomology.
\end{proof}

\begin{cor}\label{LongMayerViet}
Let $X$ be a $k$-scheme of finite type, $X_1,\ldots,X_r \subset X$ 
closed subschemes such that $X = X_1 \cup \cdots \cup X_r$. For each 
sequence $1 \leq i_0 < \cdots < i_n \leq r$, let $X_{i_0,\ldots,i_n} 
= X_{i_0} \cap \cdots \cap X_{i_n}$. Then the sequence
\eq{LongMVSeq}{ 0 \to W\sO_{X,K} \to \prod_{i=1}^r W\sO_{X_i,K} \to 
\cdots \to W\sO_{X_{1,\ldots,r},K} \to 0 }
is exact.
\end{cor}

\begin{proof}
The statement being true for $r = 2$, we proceed by induction on $r$.
Let $X' = X_2 \cup \cdots \cup X_r$. Up to a shift, the complex
\eqref{LongMVSeq} is the cone of the morphism of complexes
\eq{DecLongMVSeq}{ \xymatrix@C=0.6cm@R=0.6cm{
0 \ar[r] & W\sO_{X,K} \ar[r]\ar[d] & 
\prod_{i=2}^r W\sO_{X_i,K} \ar[r]\ar[d] & \cdots \ar[r] & 
W\sO_{X_{2,\ldots,r},K} \ar[r]\ar[d] & 0\, \\
0 \ar[r] & W\sO_{X_1,K} \ar[r] & 
\prod_{i=2}^r W\sO_{X_{1,i},K} \ar[r] & \cdots \ar[r] & 
W\sO_{X_{1,\ldots,r},K} \ar[r] & 0.
} } 
Thus it suffices to prove that this morphism is a quasi-isomorphism.
The induction hypothesis implies that the canonical morphism from the
complex $0 \to W\sO_{X,K} \to W\sO_{X',K} \to 0$ (resp.\ $0 \to
W\sO_{X_1,K} \to W\sO_{X_1\cap X',K} \to 0$) to the upper line (resp.\
lower line) of \eqref{DecLongMVSeq} is a quasi-isomorphism. Therefore
it suffices to prove that the morphism of complexes
\eq{}{ \xymatrix@C=0.6cm@R=0.6cm{
0 \ar[r] & W\sO_{X,K} \ar[r]\ar[d] & 
W\sO_{X',K} \ar[r]\ar[d] & 0 \notag\\
0 \ar[r] & W\sO_{X_1,K} \ar[r] & 
W\sO_{X_1\cap X',K} \ar[r] & 0
} } 
is a quasi-isomorphism, and this follows from Proposition \ref{MayerViet}.
\end{proof}

Our next theorem is the main result which allows to use the groups 
$H^{\ast}(X,W\sI_K)$ to define a Witt vector cohomology with compact 
supports for separated $k$-schemes of finite type. Its starting 
point is Deligne's ``independence on the compactification'' result 
in the construction of the $f_!$ functor for coherent sheaves 
\cite[App., Prop.~5]{Ha1}.

\begin{thm}\label{IndComp}
Let $f : X' \to X$ be a proper morphism between two $k$-schemes of 
finite type, and let $\sI \subset \sO_X$ be a coherent ideal, $\sI' = 
\sI\sO_{X'}$, $Y = V(\sI)$, $U = X \setminus Y$, $U' = f^{-1}(U)$.

\romain If $f$ induces a finite morphism $U' \to U$, then, for all $q 
\geq 1$,
\eq{Rqnull}{R^q f_{\ast}(W\sI'_K) = 0.}

\romain If $f$ induces an isomorphism $U' \xrightarrow{\
\sim\ \,} U$, the canonical morphism
\eq{IndCompLoc}{W\sI_K \to f_{\ast}(W\sI'_K)}
is an isomorphism, and induces an isomorphism
\eq{IndCompGlob}{\R\Gamma(X,W\sI_K) \xrightarrow{\ \sim\ \,} 
\R\Gamma(X',W\sI'_K).}
\end{thm}

The proof will use the following general lemma:

\begin{lem}\label{Wrs}
Let $I \subset A$ be an ideal in a ring $A$. For all integers $n \geq
2, r, s \in \N$, define
$$ W_n^{r,s}(I) = \{ (a_0,a_1,\ldots) \in W_n(A)\ |\ a_0 \in I^r, a_i
\in I^s \mbox{\ for all } i \geq 1\}.$$
 
\romain If $s \leq pr$, the subset $W_n^{r,s}(I)$ is an ideal of
$W_n(A)$, which sits in the short exact sequence of abelian groups
\eq{SeqWrs}{ 0 \to W_{n-1}(I^s) \xrightarrow{V} 
W_n^{r,s}(I) \xrightarrow{R^{n-1}} I^r \to 0. }

\romain If $r \leq s \leq pr$, the sequence \eqref{SeqWrs} sits in a 
commutative diagram of short exact sequences
\eq{DiagWrs}{\xymatrix{
0\ar[r] & W_{n-1}(I^s)\ar@{=}[d]\ar[r]^-V & 
W_n(I^s)\ar@{^{(}->}[d]\ar[r]^-{R^{n-1}} & 
I^s\ar@{^{(}->}[d]\ar[r] & 0\\
0\ar[r] & W_{n-1}(I^s)\ar@{^{(}->}[d]\ar[r]^-V & 
W_n^{r,s}(I)\ar@{^{(}->}[d]\ar[r]^-{R^{n-1}} & 
I^r\ar@{=}[d]\ar[r] & 0\\
0\ar[r] & W_{n-1}(I^r)\ar[r]^-V & 
W_n(I^r)\ar[r]^-{R^{n-1}} & I^r\ar[r] & 0
}}
where the vertical arrows are the natural inclusions.
\end{lem}

\begin{proof}
Assume first that $s \leq pr$. The subset $W_n^{r,s}(I)$ can be
described as the set of Witt vectors of the form $\ua + V(b)$, where
$a \in I^r$, $b \in W_{n-1}(I^s)$, and $\underline{a}$ denotes the 
Teichm\"uller representative of $a$. To prove that it is an
additive subgroup, it suffices to verify that if $a, a' \in I^r$,
then $\ua + \ua' = \underline{a+a'}+V(c)$, with $c \in W_{n-1}(I^s)$.
If $S_i(X_0,\ldots,X_i,Y_0,\ldots,Y_i)$ are the universal polynomials
defining the addition in $W_n(A)$, it is easy to check that
$S_i(X_0,0,\ldots,Y_0,0,\ldots)$ is a homogeneous polynomial of
degree $p^i$ in $\Z[X_0,Y_0]$ \cite[II \S 6]{Se3}. Since $pr \geq s$,
the claim follows. Using again the condition $pr \geq s$, the fact
that $W_n^{r,s}(I)$ satisfies the multiplicativity property of an
ideal follows from the relations
$$ \ua(x_0,x_1,x_2,\ldots) = (ax_0,a^px_1,a^{p^2}x_2,\ldots), \quad\quad
V(b)x = V(bF(x)). $$

If we assume in addition that $r \leq s$, then the vertical 
inclusions of the diagram are defined, and its commutativity is 
obvious.
\end{proof}

We will use repeatedly the following elementary remark.

\begin{lem}\label{Butterfly}
Let $\sC$ be an abelian category, and let
\eq{DiagGn}{ \xymatrix{
& E\ar[d]^{v}\ar[rd]^{u''} & \\
G'\ar[dr]_{u'}\ar[r]^-{V} & 
G\ar[d]^{w}\ar[r]^-{R} & G''\\
& F & 
} }
be a commutative diagram of morphisms of $\sC$ such that the horizontal 
sequence is exact. If $u'=0$ and $u''=0$, then $w \circ v = 0$.
\end{lem}

\begin{proof}
Exercise.
\end{proof}

\begin{lem}\label{EssZero}\setcounter{romain}{1}
\textnormal{(i)}\ \ Under the assumptions of Theorem \ref{IndComp},  (i)
there exists an integer $a \geq 0$ such that, for all $q \geq 1$, all
$n \geq 1$ and all $r \geq 0$, the canonical morphism
\eq{ZeroMor}{R^qf_{\ast}(W_n(\sI'{}^{r+a})) \to R^qf_{\ast}(W_n(\sI'{}^r))}
is the zero morphism.

\romain Assume that $f$ induces an isomorphism $U' \xrightarrow{\
\sim\ \,} U$, and define
\ga{}{\sK_n^r = \Ker\big(W_n(\sI^r) \to f_{\ast}(W_n(\sI'{}^r))\big), \notag\\
\sC_n^r = \Coker\big(W_n(\sI^r) \to f_{\ast}(W_n(\sI'{}^r))\big). \notag}
Then there exists an integer $a \geq 0$ such that, for all
$n \geq 1$ and all $r \geq 0$, the canonical morphisms
\eq{ZeroMor0}{\sK_n^{r+a} \to \sK_n^r, \quad\quad \sC_n^{r+a} \to
\sC_n^r,}
are the zero morphisms.
\end{lem}
\begin{proof}
We first prove assertion (i). We may fix $q \geq 1$, since
$R^qf_{\ast} = 0$ for $q$ big enough. 
When $n = 1$, we can apply Deligne's result \cite[App.,
Prop.~5]{Ha1} to $\sO_{X'}$, and we obtain an integer $b \geq 0$
such that, for all $r \geq 0$, the canonical morphism
$$ u_1 : R^q f_{\ast}(\sI'{}^{r+b}) \to R^q f_{\ast}(\sI'{}^r) $$
is $0$ (note that, for this result, Deligne's argument only uses that
$f$ is finite above $U$). Let us prove by induction on $n$ that,
for all $n \geq 1$ and all $r \geq b$, the canonical morphism
$$ u_n : R^q f_{\ast}(W_n(\sI'{}^{r+b})) \to R^q f_{\ast}(W_n(\sI'{}^r)) $$
is also $0$. 

The condition $r \geq b$ implies that the couple $(r,r+b)$ is such
that $r \leq r+b \leq pr$. Therefore, we may use Lemma \ref{Wrs} to
define ideals $W_n^{r,r+b}(\sI') \subset W_n(\sO_{X'})$, and, for $n
\geq 2$, we obtain a commutative diagram \eqref{DiagWrs} relative to
$\sI'$ and $(r, r+b)$. Since the middle row of \eqref{DiagWrs} is a 
short exact sequence, the diagram
$$ \xymatrix{
& 
R^q f_{\ast}(W_n(\sI'{}^{r+b}))\ar[d]\ar[r]^-{R^{n-1}} & 
R^q f_{\ast}(\sI'{}^{r+b})\ar[d]^{u_1}\\
R^q f_{\ast}(W_{n-1}(\sI'{}^{r+b}))\ar[d]^{u_{n-1}}\ar[r]^-V & 
R^q f_{\ast}(W_n^{r,r+b}(\sI'))\ar[d]\ar[r]^-{R^{n-1}} & 
R^q f_{\ast}(\sI'{}^r)\\
R^q f_{\ast}(W_{n-1}(\sI'{}^r))\ar[r]^-V & 
R^q f_{\ast}(W_n(\sI'{}^r)) & 
} $$ 
obtained by applying the functor $R^q f_{\ast}$ to \eqref{DiagWrs} has
an exact middle row. As $u_1 = 0$, and the composition of the middle
vertical arrows is $u_n$, Lemma \ref{Butterfly} shows by induction that
$u_n = 0$ for all $n \geq 1$. If we set $a = 2b$, then assertion (i)
holds.

Assertion (ii) can be proved by the same method. Let $\sR =
\bigoplus_{r \geq 0} \sI^r$. It follows from \cite[3.3.1]{EGA} that
the quasi-coherent graded $\sR$-modules $\sK_1 = \bigoplus_{r \geq 0}
\sK_1^r$ and $\sC_1 = \bigoplus_{r \geq 0} \sC_1^r$ are finitely
generated. As the restriction of $f$ to $U'$ is an isomorphism, they
are supported in $Y$. Therefore, there exists an integer $m$ such
that $\sI^{m} \sK_1^r = \sI^{m} \sC_1^r = 0$ for all $r \geq 0$.
Moreover, there exists an integer $d$ such that, for all $r \geq d$
the canonical morphisms
$$ \sI \otimes_{\sO_X} \sK_1^r \to \sK_1^{r+1}, \quad 
\sI \otimes_{\sO_X} \sC_1^r \to \sC_1^{r+1}, $$
are surjective. It follows that the morphisms
$$ \sK_1^{r+m} \to \sK_1^r, \quad  \sC_1^{r+m} \to \sC_1^r $$ 
are $0$ for $r \geq d$. Replacing $m$ by $b = d + m$, the
corresponding morphisms are $0$ for all $r \geq 0$, which proves
assertion (ii) when $n = 1$.

As $\sK_1^r \subset \sO_X$ for all $r$, this implies that $\sK_1^r = 0$ 
for $r \geq b$. Thanks to the exact sequences 
$$ 0 \to \sK_{n-1}^r \xrightarrow{V} \sK_n^r \xrightarrow{R^{n-1}} 
\sK_1^r, $$
it follows that $\sK_n^r=0$ for all $n \geq 1$ and all $r \geq b$. 
Thus assertion (ii) holds for the modules $\sK_n^r$, with $a = b$. 

To prove it for the modules $\sC_n^r$, we introduce for $r \geq b$ 
the modules 
$$ \sC_n^{r,r+b} = \Coker\big(W_n^{r,r+b}(\sI) \to 
f_{\ast}(W_n^{r,r+b}(\sI'))\big). $$
The snake lemma applied to the diagrams
$$ \xymatrix{
0 \ar[r] & W_{n-1}(\sI^{r+b}) \ar[r]^-V \ar[d] & 
W_n^{r,r+b}(\sI) \ar[r]^-{R^{n-1}} \ar[d] & \sI^r \ar[r] \ar[d] & 0 \\
0 \ar[r] & f_{\ast}(W_{n-1}(\sI'{}^{r+b})) \ar[r]^-V & 
f_{\ast}(W_n^{r,r+b}(\sI')) \ar[r]^-{R^{n-1}} & 
f_{\ast}(\sI'{}^r) 
} $$ 
gives exact sequences
$$ \sC_{n-1}^{r+b} \xrightarrow{V} \sC_n^{r,r+b} 
\xrightarrow{R^{n-1}} \sC_1^r. $$ 
The functoriality of the diagrams \eqref{DiagWrs} imply that these 
exact sequences sit in commutative diagrams
$$ \xymatrix{
& \sC_n^{r+b}\ar[d]\ar[r]^-{R^{n-1}} & 
\sC_1^{r+b}\ar[d]^{u_1}\\
\sC_{n-1}^{r+b}\ar[d]^{u_{n-1}}\ar[r]^-V & 
\sC_n^{r,r+b}\ar[d]\ar[r]^-{R^{n-1}} & 
\sC_1^r\\
\sC_{n-1}^r\ar[r]^-V & 
\sC_n^r & \quad,
} $$ 
where the morphisms $u_i$ are the canonical morphisms, and the 
composition of the middle vertical arrows is $u_n$. As $u_1 = 0$, Lemma 
\ref{Butterfly} implies by induction that $u_n = 0$ for all $n$. This 
proves assertion (ii) for $\sC_n^r$, with $a = 2b$.
\end{proof}

\subsection{}\label{ProofInd} \textit{Proof of Theorem \ref{IndComp}.}
Under the assumptions of \ref{IndComp}, let $q \geq 1$ be an integer,
and let $a$ be an integer satisfying the conclusion of Lemma
\ref{EssZero} (i) for the family of sheaves $R^q
f_{\ast}(W_n(\sI'{}^r))$, for all $n \geq 1$ and $r \geq 0$. Let $c$
be such that $p^c > a$. Since the Frobenius map $F^c : R^q
f_{\ast}(W_n(\sI'{}^r)) \to R^q f_{\ast}(W_n(\sI'{}^r))$ factors
through $R^q f_{\ast}(W_n(\sI'{}^{p^c r}))$, it follows from Lemma
\ref{EssZero} that $F^c$ acts by zero on $R^q
f_{\ast}(W_n(\sI'{}^r))$, for all $n \geq 1$ and all $r \geq 1$.
Therefore, for all $r \geq 1$, $F^c$ acts by zero on $\R\varprojlim_n
(R^q f_{\ast}(W_n(\sI'{}^r)))$.

In particular, $F^c$ acts by zero on the sheaves $\varprojlim_n
R^qf_{\ast}(W_n\sI')$ and $R^1\varprojlim_n R^qf_{\ast}(W_n\sI')$. As
the inverse system $(W_n\sI')_{n \geq 1}$ has surjective transition
maps, and terms with vanishing cohomology on affine open subsets, it
is $\varprojlim_n$-acyclic, and we obtain
$$ \R f_{\ast}(W\sI') = \R f_{\ast}(\R\varprojlim_n W_n\sI') 
\cong \R\varprojlim_n \R f_{\ast}(W_n\sI'). $$
This isomorphism provides a biregular spectral sequence 
$$ E_2^{i,j} = R^i\varprojlim_n R^j f_{\ast}(W_n\sI') \Rightarrow 
R^{i+j}f_{\ast}(W\sI'), $$
in which the filtration of the terms $R^{i+j}f_{\ast}(W\sI')$ is
of length $\leq 2$ since the functors $R^i\varprojlim_n$ are
zero for $i \geq 2$. As $F^c$ acts by zero on the terms $E_2^{i,j}$ 
for $j \geq 1$, $F^{2c}$ acts by zero on the terms $R^{i+j}f_{\ast}(W\sI')$ 
for $i+j \geq 2$. 

On the other hand, the term $E_2^{1,0} = R^1\varprojlim_n
f_{\ast}(W_n\sI')$ is $0$, because the morphisms
$f_{\ast}(W_{n+1}\sI') \to f_{\ast}(W_n\sI')$ are surjective (since
$f_{\ast}(W_n\sI') \cong W_n(f_{\ast}(\sI'))$), and the cohomology of
the terms $f_{\ast}(W_n\sI')$ vanishes on any open affine subset.
Therefore, $F^c$ acts by zero on $R^1f_{\ast}(W\sI')$.

Thus the action of $F$ on $R^q f_{\ast}(W\sI')$ is nilpotent for all 
$q \geq 1$. Since this action becomes an 
isomorphism after tensorisation with $K$, we obtain that $R^q 
f_{\ast}(W\sI'_K) = R^q f_{\ast}(W\sI')_K = 0$ for all $q \geq 1$, 
which proves assertion (i) of Theorem \ref{IndComp}. 

Let us assume now that $f : U' \to U$ is an isomorphism. It follows
from Lemma \ref{EssZero} (ii) that there exists an integer $a$ such
that, for all $n \geq 1$ and all $r \geq 0$, the morphisms
$\sK_n^{r+a} \to \sK_n^r$ and $\sC_n^{r+a} \to \sC_n^r$ are zero.
Taking $c$ such that $p^c > a$, it follows that $F^c$ acts by zero on
$\sK_n^r$ and $\sC_n^r$ for all $n$ and all $r$. Therefore, $F^c$ acts
by zero on $\varprojlim_n \sK_n^1$ and $\varprojlim_n \sC_n^1$. It is
easy to check that these two inverse systems and the inverse system
$\Im\big(W_n\sI \to f_{\ast}(W_n\sI')\big)$ all have surjective
transition maps. As their terms have vanishing cohomology on affine
open subsets, they are $\varprojlim$-acyclic, and we obtain
\ga{}{\Ker\big(W\sI \to f_{\ast}(W\sI')\big) =  \varprojlim_n 
\sK_n^1, \notag\\
\Coker\big(W\sI \to f_{\ast}(W\sI')\big) =  \varprojlim_n 
\sC_n^1. \notag}
After tensorisation with $K$, $F$ becomes an isomorphism on
$\varprojlim_n \sK_n^1$ and $\varprojlim_n \sC_n^1$, and assertion 
(ii) of Theorem \ref{IndComp} follows. \hfill $\Box$

\subsection{}\label{CohCompSupp}
We now observe that the previous results imply that the cohomology
spaces $H^{\ast}(X, W\sI_K)$ only depend on the $k$-scheme $U$, and
have the same functoriality properties with respect to $U$ than
cohomology with compact supports.

Indeed, suppose $U$ is fixed, and let $U \inj X_1, U
\inj X_2$ be two open immersions into proper $k$-schemes.
Let $X \subset X_1 \times_k X_2$ be the scheme theoretic closure of
$U$ embedded diagonally into $X_1 \times_k X_2$. The two projections
induce proper maps $p_1 : X \to X_1, p_2 : X \to X_2$. If $\sI \subset
\sO_X$, $\sI_1 \subset \sO_{X_1}$, $\sI_2 \subset \sO_{X_2}$ are the
ideals defining $Y = (X \setminus U)_{\red}$, $Y_1 = (X_1 \setminus
U)_{\red}$ and $Y_2 = (X_2 \setminus U)_{\red}$, we deduce from
\ref{TopInv} and \ref{IndComp} that the homomorphisms
$$ \R\Gamma(X_2, W\sI_{2,K}) \xrightarrow{p_2^{\ast}} \R\Gamma(X, 
W\sI_K) \xleftarrow{p_1^{\ast}} \R\Gamma(X_1, W\sI_{1,K}) $$ 
are isomorphisms. If we define 
$$ \varepsilon_{12} = p_1^{\ast\;-1} \circ p_2^{\ast} : \R\Gamma(X_2,
W\sI_{2,K}) \xrightarrow{\sim} \R\Gamma(X_1, W\sI_{1,K}), $$
it is easy to check that the isomorphisms $\varepsilon_{ij}$ satisfy
the transitivity condition for a third open immersion $U
\inj X_3$ into a proper $k$-scheme. Therefore, they provide
canonical identifications between the cohomology complexes
$\R\Gamma(X, W\sI_K)$ defined by various open immersions of $U$ into
proper $k$-schemes $X$.

Assume $U$ is a separated $k$-scheme of finite type. Since, by
Nagata's theorem, there exists a proper $k$-scheme $X$ and an open
immersion $U \inj X$, this independence property allows to
define the \textit{Witt vector cohomology with compact supports} of
$U$ by setting 
\ga{defHcW}{ \R\Gamma\cpt(U, W\sO_{U,K}) := \R\Gamma(X,W\sI_K) 
\cong \R\Gamma(X,W\sI)_K, \\
H^{\ast}\cpt(U, W\sO_{U,K}) := H^{\ast}(X, W\sI_K) 
\cong H^{\ast}(X, W\sI)_K, }
where $\sI \subset \sO_X$ is any coherent ideal defining the closed
subset $X \setminus U$. As the restriction of $W\sI$ to $U$ is $W\sO_U$, 
there is a canonical morphism $\R\Gamma\cpt(U, W\sO_{U,K}) \to 
\R\Gamma(U, W\sO_{U,K})$, which is an isomorphism when $U$ is proper.

These cohomology groups have the following functoriality properties:

\romain They are contravariant with respect to proper maps. 

Let $f : U' \to U$ be a proper $k$-morphism of separated $k$-schemes,
and let $U' \inj X', U \inj X$ be open
immersions into proper $k$-schemes. Replacing if necessary $X'$ by the
scheme theoretic closure of the graph of $f$ in $X' \times_k X$, we
may assume that there exists a $k$-morphism $g : X' \to X$ extending
$f$, and that $U'$ is dense in $X'$. As $f$ is proper, it follows that
$U'= g^{-1}(U)$. If $\sI \subset \sO_X$ is any coherent ideal such 
that $V(\sI) = X \setminus U$, then $\sI' = \sI\sO_{X'}$ is such that 
$V(\sI') = X' \setminus U'$, and we can define the homomorphism
$$ f^{\ast} : \R\Gamma\cpt(U, W\sO_{U,K}) \to 
\R\Gamma\cpt(U', W\sO_{U',K})$$
as being 
$$ g^{\ast} : \R\Gamma(X, W\sI_K) \to \R\Gamma(X', W\sI'_K).$$
We leave as an exercise to check that, up to canonical isomorphism, 
$f^{\ast}$ does not depend on the choices.

\romain They are covariant with respect to open immersions. 

Let $j : V \inj U$ be an open immersion, let
$U \inj X$ be an open immersion into a proper $k$-scheme, and let
$\sI, \sJ \subset \sO_X$ be the ideals of $Y = (X \setminus U)_{\red}$ and
$Z = (X \setminus V)_{\red}$. Then the homomorphism
$$ j_{\ast} : \R\Gamma\cpt(V, W\sO_{V,K}) \to 
\R\Gamma\cpt(U, W\sO_{U,K}) $$
is defined as being
$$ \R\Gamma(X, W\sJ_K) \to \R\Gamma(X, W\sI_K). $$
If $T = U \setminus V = Z \cap U$, then $T$ is open in the proper
$k$-scheme $Z$, and the ideal $\sI/\sJ \subset \sO_Z$ defines the
complement of $T$ in $Z$. Thus, the usual distinguished triangle
\eq{TrHcW}{ \R\Gamma\cpt(V, W\sO_{V,K}) \to 
\R\Gamma\cpt(U, W\sO_{U,K}) \to 
\R\Gamma\cpt(T, W\sO_{T,K}) \xrightarrow{+1} }
is obtained by tensoring the short exact sequence
$$ 0 \to W\sJ \to W\sI \to W(\sI/\sJ) \to 0, $$
with $K$, and taking its cohomology on $X$.

\begin{prop}\label{finite}
For any separated $k$-scheme of finite type $U$, the cohomology spaces 
$H^{\ast}\cpt(U, W\sO_{U,K})$ are finite dimensional $K$-vector 
spaces, on which the Frobenius endomorphism has slopes in $[0,1[$.
\end{prop}

\begin{proof}
It suffices to prove the statement when $U$ is a proper $k$-scheme
$X$. Then the statement reduces to the finiteness of the usual
cohomology spaces $H^{\ast}(X, W\sO_X)_K$. This is a well-known result
(\cf \cite[III, Th. 2.2]{B}, which is valid for $W\sO_X$ without the
smoothness assumption). For the sake of completeness, we give a proof
here.

Write $M:= H^i(X, W\sO_X)$. As $X$ is proper, the groups $H^i(X,
W_n\sO_X)$ satisfy the Mittag-Leffler condition, so that the
homomorphism
$$ H^i(X, W\sO_X) \to \varprojlim_n H^i(X, W_n\sO_X) $$ 
is an isomorphism. Therefore, the Verschiebung endomorphism endows $M$
with a structure of module over the non necessarily commutative ring
$R:= W_{\sigma}[[V]]$, where the index $\sigma$ refers to the
commutation rules $aV=V\sigma(a)$ for $a \in W$. Then $M$ is separated
and complete for the $V$-adic topology. On the other hand, $M/VM \inj
H^i(X,\sO_X)$ is a finite dimensional $k$-vector space. Since $M/VM$
is finitely generated over $W$, it follows that $M$ is finitely
generated over $R$. Moreover, the fact that $M/VM$ has finite length
implies that $M$ is a torsion $R$-module. Thus, there exists a finite
number of non zero elements $a_i(V) = a_{i,0} + a_{i,1}V + a_{i,2}V^2
+ \ldots \in R, i = 1,\ldots,r,$ and a surjection
\eq{PresM}{\bigoplus_{i=1}^r R/Ra_i(V) \surj M.}
Fix some $i, 1 \leq i \leq r$. We are interested in the module
structure after inverting $p$, so we may assume that $a_{i,\ell} \in
W(k)^\times$ for some $\ell$. Let $\ell$ be minimal, so $p|a_{i, j},\
0\le j <\ell$. Even if the ring $R$ is not commutative, we can
get a factorization
$$ \sum_j a_{i,j}V^j = \Big(c_{i,0}+c_{i,1}V+\ldots\Big)\Big(b_{i,0}+
\ldots + b_{i,\ell-1}V^{\ell-1} + V^\ell\Big); \ \ p|b_{i,j}. $$
To see this, we factor
$$ \sum_j a_{i,j}V^j \equiv
\Big(c_{i,0}^{(s)}+c_{i,1}^{(s)}V+\ldots\Big)\Big(b_{i,0}^{(s)}+\ldots +
b_{i,\ell-1}^{(s)}V^{\ell-1} + V^\ell\Big)\ \ \mathrm{mod} \ p^{s+1}, $$
where $b_{i,j}^{(s+1)}, c_{i,j}^{(s+1)} \equiv b_{i,j}^{(s)},
c_{i,j}^{(s)} \mod p^{s+1}$, starting with $b_{i,j}^{(0)} = 0$ and
$c_{i,j}^{(0)} = a_{i,j+\ell}$. Note that $c_{i,0}^{(s)} \in
W(k)^\times$. Details are standard and are left for the reader.
Writing $b_i(V) = b_{i,0} + b_{i,1}V + \ldots + V^\ell$ with $b_{i,j}
= \lim b_{i,j}^{(s)}$, it follows that $R/Ra_i(V) \cong R/Rb_i(V)
\cong W(k)^{\oplus \ell}$. The finiteness now follows from \eqref{PresM}.

As $H^{\ast}(X, W\sO_X)$ is a finitely generated $W$-module modulo
torsion, the slopes of the Frobenius endomorphism are positive. On 
the other hand, the existence of a Verschiebung endomorphism $V$ 
such that $FV = VF = p$ implies that all slopes are $\leq 1$. As $V$ 
is topologically nilpotent, there cannot be any non zero element of 
slope $1$ in $H^{\ast}(X, W\sO_{X,K})$, and all its slopes belong to 
$[0,1[$. 
\end{proof}

\bigskip
\section{A descent theorem}
\label{SectDesc}
\medskip

We will need simplicial resolutions based on de Jong's fundamental
result on alterations. In this section, we briefly recall some
related definitions, and how to construct such resolutions. We then
prove for Witt vector cohomology with compact supports a particular 
case of \'etale cohomological descent which will be one of the 
main ingredients in our proof of Theorem \ref{SlopeHc}.

If $X$ is a scheme, and $n \geq -1$ an integer, we denote
as usual by $\sk^X_n$ the truncation functor from the category of
simplicial $X$-schemes to the category of $n$-truncated simplicial
$X$-schemes, and by $\cosk^X_n$ its right adjoint \cite[5.1]{De1}.

\begin{defn}\label{resolutions} Let $X$ be a reduced $k$-scheme of
finite type, and let $X\lbul$ be a simplicial $k$-scheme (resp.\
$N$-truncated simplicial scheme, for some $N \in \N$).

\alphab A $k$-augmentation $f\lbul : X\lbul \to X$ is called a
\textit{proper hypercovering} (resp.\ \textit{$N$-truncated
proper hypercovering}) of $X$ if, for all $n \geq 0$ (resp. $0 \leq n
\leq N$), the canonical morphism $X_n \to
\cosk^X_{n-1}(\sk^X_{n-1}(X\lbul))_n$ is proper and surjective.

\alphab A $k$-augmentation $f\lbul : X\lbul \to X$ is called an
\textit{\'etale hypercovering} (resp.\ \textit{$N$-truncated
\'etale hypercovering}) of $X$ if, for all $n \geq 0$ (resp. $0 \leq n
\leq N$), the canonical morphism $X_n \to
\cosk^X_{n-1}(\sk^X_{n-1}(X\lbul))_n$ is \'etale and surjective.

\alphab A $k$-augmentation $f\lbul : X\lbul \to X$ is
called a \textit{simplicial resolution} (resp.\ \textit{$N$-truncated
simplicial resolution}) of $X$ if the following conditions hold:

\romain $f\lbul$ is a proper hypercovering (resp.\ $N$-truncated
proper hypercovering) of $X$.

\romain For all $n \geq 0$ (resp.\ $0 \leq n \leq N$), there exists a
dense open subset $U_n \subset X$ such that the restriction of
$\sk^X_n(X\lbul)$ above $U_n$ is a $n$-truncated \'etale hypercovering
of $U_n$.

\romain For all $n \geq 0$ (resp. $0 \leq n \leq N$), $X_n$ is a smooth 
quasi-projective $k$-scheme.
\end{defn} 

Since $k$ is perfect, de Jong's theorem \cite[Th.\ 4.1 and Rem.\
4.2]{dJ} implies that, for any separated integral (\ie reduced and
irreducible) $k$-scheme of finite type $X$, there exists a
quasi-projective, smooth and integral $k$-scheme $X'$, and a
surjective, projective and generically \'etale morphism $X' \to X$. If
$X$ is separated and reduced, but not necessarily irreducible, one can
apply de Jong's theorem to each irreducible component of $X$, and one
gets in this way a quasi-projective smooth $k$-scheme $X'$, and a
surjective projective morphism $X' \to X$ which is \'etale over a
dense open subset of $X$. Replacing resolution of singularities by
this result, one can then proceed as Deligne \cite[(6.2.5)]{De1} to
show that any separated reduced $k$-scheme of finite type has a
simplicial resolution.
\medskip

\noindent\textit{Remark.} Let $\sC$ be the category of $X$-schemes
which are separated of finite type over $k$, and $\sD$ the subcategory
defined as follows:
\begin{list}{-}{\setlength{\leftmargin}{6mm}\setlength{\topsep}{0mm}}
\item $\Ob(\sD) = \Ob(\sC)$;
\item If $Y, Z \in \Ob(\sD)$, $\Hom_{\sD}(Y, Z)$ is the set of
$X$-morphisms $f : Y \to Z$ which are proper, surjective, and such
that there exists a dense open subset $U \subset X$ with the property
that the restriction $f|_U : Y|_U \to Z|_U$ of $f$ above $U$ is
\'etale.
\end{list}
\noindent Then $\sD$ satisfies the condition (HC) of \cite[5.1]{Ts2},
and the simplicial resolutions of $X$ are the $\sD$-hypercoverings of
$X$ by smooth quasi-projective $k$-schemes, in the sense of
\cite[5.1.1]{Ts2}.
\medskip

\begin{prop}\label{descent}
Let $X$ be a reduced $k$-scheme of finite type, let $N \in \N$, and let
$f\lbul : X\lbul \to X$ be a $N$-truncated proper hypercovering of
$X$. Let $\sI \subset \sO_X$ be a coherent ideal, $\sI_q =
\sI\sO_{X_q}$, $Y = V(\sI)$, $U = X \setminus Y$, $U_q = f_q^{-1}(U)$.
Assume that the restriction $U\lbul$ of $X\lbul$ above $U$ is an 
$N$-truncated \'etale hypercovering of $U$. Then:

\romain The complex
$$ 0 \to W\sI_K \to f_{0\,\ast}(W\sI_{0,K}) \to 
f_{1\,\ast}(W\sI_{1,K}) \to 
\cdots \to  f_{N\,\ast}(W\sI_{N,K}), $$
where $W\sI_K$ sits in degree $-1$, is acyclic in degrees $\neq N$.

\romain The canonical morphisms 
\eq{CohDesc}{ H^q(X, W\sI_K) \to H^q(X\lbul, W\sI\lbul{}_{\!,K}) }
are isomorphisms for all $q < N$.
\end{prop}

\begin{proof}\setcounter{romain}{0}
For each $n \geq 1$ and each $r \geq 0$, we denote by 
$\sL_n^{\bul,r}$ the complex
\eq{DefsLr}{ 0 \to W_n(\sI^r) \to f_{0\,\ast}(W_n(\sI_0^r)) \to 
\cdots \to  f_{N\,\ast}(W_n(\sI_N^r)), }
where $W_n(\sI^r)$ sits in degree $-1$, and we set
\ga{DefsL}{ \sL_n\hbul = \bigoplus_{r \geq 0} \sL_n^{\bul,r}, \\
\label{DefsH}\sH_n^{q,r} = \sH^q(\sL_n^{\bul,r}), \quad 
\sH^q_n = \sH^q(\sL_n\hbul) = \bigoplus_{r \geq 0} \sH_n^{q,r}. } 

If we denote again by $\sR$ the graded $\sO_X$-algebra $\bigoplus_{r
\geq 0} \sI^r$, it follows from \cite[3.3.1]{EGA} that $\sL_1\hbul$ is
a complex of quasi-coherent graded $\sR$-modules of finite type.
Therefore, $\sH^q_1$ is a quasi-coherent graded $\sR$-module of finite
type for all $q$, and there exists an integer $d \geq 0$ such that the
morphism
$$ \sI \otimes_{\sO_X} \sH_1^{q,r} \to \sH_1^{q,r+1} $$
is surjective for all $r \geq d$. Moreover, $U\lbul$ is an
$N$-truncated \'etale hypercovering of $U$, and therefore it satisfies
cohomological descent for quasi-coherent modules. It follows that
$\sL_1\hbul|_U$ is acyclic in degrees $< N$. Then we can find an
integer $m \geq 0$ such that $\sI^m\sH_1^{q,r} = 0$ for all $q < N$
and all $r \geq 0$. For $r \geq d$, this implies that the image of
$\sH_1^{q,r+m}$ in $\sH_1^{q,r}$ is $0$. Finally, if we set $b = d+m$,
we obtain that the canonical morphism $\sH_1^{q,r+b} \to \sH_1^{q,r}$
is $0$ for all $q < N$ and all $r \geq 0$.

We can now proceed as in the proof of Lemma \ref{EssZero} to prove
that the morphism $\sH_n^{q,r+b} \to \sH_n^{q,r}$ is $0$ for all $q
\leq N$, all $n \geq 1$ and all $r \geq b$. As this last condition
implies that $r \leq r+b \leq pr$, we can introduce for $n \geq 2$ the
subcomplex $\sL_n^{\bul,r,r+b} \subset \sL_n^{\bul,r}$ defined by
$$ 0 \to W_n^{r,r+b}(\sI) \to f_{0\,\ast}(W_n^{r,r+b}(\sI_0)) \to 
\cdots \to f_{N\,\ast}(W_n^{r,r+b}(\sI_N)), $$
and $\sH_n^{q,r,r+b} = \sH^q(\sL_n^{\bul,r,r+b})$. From Lemma 
\ref{Wrs}, we deduce a commutative diagram of morphisms of complexes 
$$ \xymatrix{
0\ar[r] & \sL_{n-1}^{\bul,r+b}\ar@{=}[d]\ar[r]^-V & 
\sL_n^{\bul,r+b}\ar@{^{(}->}[d]\ar[r]^-{R^{n-1}} & 
\sL_1^{\bul,r+b}\ar@{^{(}->}[d]\ar[r] & 0\,\\
0\ar[r] & \sL_{n-1}^{\bul,r+b}\ar@{^{(}->}[d]\ar[r]^-V & 
\sL_n^{\bul,r,r+b}(I)\ar@{^{(}->}[d]\ar[r]^-{R^{n-1}} & 
\sL_1^{\bul,r}\ar@{=}[d]\ar[r] & 0\,\\
0\ar[r] & \sL_{n-1}^{\bul,r}\ar[r]^-V & 
\sL_n^{\bul,r}\ar[r]^-{R^{n-1}} & \sL_1^{\bul,r}\ar[r] & 0,
} $$
where the rows are short exact sequences, because any section of a 
sheaf $f_{q\,\ast}(\sI_q^r)$ can be lifted to a section of 
$f_{q\,\ast}(W_n(\sI_q^r))$ or $f_{q\,\ast}(W_n^{r,r+b}(\sI_q))$ by 
taking its Teichm\"uller representative. From the corresponding 
diagram of cohomology exact sequences, we can 
extract for each $q$ the commutative diagram
$$ \xymatrix{
& \sH_n^{q,r+b}\ar[d]\ar[r]^-{R^{n-1}} & 
\sH_1^{q,r+b}\ar[d]^{u_1}\\
\sH_{n-1}^{q,r+b}\ar[d]^{u_{n-1}}\ar[r]^-V & 
\sH_n^{q,r,r+b}\ar[d]\ar[r]^-{R^{n-1}} & 
\sH_1^{q,r}\\
\sH_{n-1}^{q,r}\ar[r]^-V & 
\sH_n^{q,r} & \quad\quad\ \ ,
} $$
where the morphisms $u_1$, $u_{n-1}$ and the composition of the two
vertical morphisms are the canonical morphisms, and the midle row is
exact. Since $u_1 = 0$ when $q < N$, Lemma \ref{Butterfly} shows by
induction that $u_n = 0$ for all $q < N$, all $n\geq 1$ and all $r
\geq b$.

Setting $a = 2b$, we obtain that the morphisms $\sH_n^{q,r+a} \to
\sH_n^{q,r}$ are $0$ for all $q < N$, all $n \geq 1$ and all $r \geq
0$. This implies that, if $c$ is chosen so that $p^c > a$, then $F^c$
acts by $0$ on the sheaves $\sH_n^{q,r}$ for $r \geq 1$. Therefore, 
$F^c$ acts by $0$ on $R^i\varprojlim_n\sH_n^{q,r}$ for all $i \geq 0$, all 
$q < N$ and all $r \geq 1$. 

Let $\sL\hbul$ be the complex
$$ 0 \to W\sI \to f_{0\,\ast}(W\sI_0) \to \cdots \to 
f_{N\,\ast}(W\sI_N).
$$ 
As the sheaves $\sL_n^{q,1}$ have vanishing cohomology on affine open
subsets of $X$, and the transition maps $\sL_{n+1}^{q,1} \to
\sL_n^{q,1}$ are surjective, the projective systems $(\sL_n^{q,1})_{n
\geq 1}$ are $\varprojlim_n$-acyclic, and the canonical morphism
$$ 
\sL\hbul \to \R\varprojlim_n \sL_n^{\bul,1}
$$ 
is an isomorphism. This isomorphism provides a spectral sequence
$$ 
E_2^{i,j} = R^i\varprojlim_n \sH^j(\sL_n\hbul{}^{\!,1}) \Rightarrow 
\sH^{i+j}(\sL\hbul),
$$ 
on which $F$ acts. As $F$ is nilpotent on the terms $E_2^{i,j}$ for $j 
< N$, it follows that $F$ is also nilpotent on $\sH^q(\sL\hbul)$ for 
$q < N$. But $F$ becomes an isomorphism after tensorisation by $K$, 
so this implies assertion (i).

Thanks to Theorem \ref{IndComp}, $R^i f_{q\,\ast}(W\sI_q)_K = 0$ for 
all $q$ and all $i \geq 1$. Therefore, the canonical morphism
$$ f\lbul{}_{\ast}(W\sI\lbul)_K \to \R f\lbul{}_{\ast}(W\sI\lbul)_K $$
is an isomorphism. Then (ii) follows from (i).
\end{proof}

\bigskip
\section{Witt vector cohomology and rigid cohomology}
\label{RigandWitt}
\medskip

In this section, we first recall how to compute rigid cohomology in
terms of de Rham complexes for the Zariski topology. Using this
description, we construct the canonical morphism from rigid cohomology
with compact supports to Witt vector cohomology with compact supports.

\subsection{}\label{tubes}  
We begin by recalling the general construction of tubes \cite{Be0},
which is sufficient to define the rigid cohomology groups of a proper
$k$-scheme.

Let $X$ be a separated $k$-scheme of finite type, and let $X \inj \P$
be a closed immersion into a smooth formal $W$-scheme. The formal
scheme $\P$ has a generic fibre $\P_K$, which is a rigid analytic
space, endowed with a continuous morphism $\spm : \P_K \to \P$, the
specialization morphism, such that $\spm^{-1}(\U) = \U_K$ for any open
subset $\U \subset \P$, and $\spm_{\ast}\sO_{\P_K} = \sO_{\P,K}$. When
$\P$ is an affine formal scheme $\Spf A$, where $A$ is a $W$-algebra
which is topologically of finite type, its generic fibre is the
affinoid space $\Spm A_K$, defined by the Tate algebra $A_K$.

If $\P = \Spf A$, and if $f_1,\ldots,f_r \in A$ is a
family of generators of the ideal $\sJ$ of $X$ in $\P$, the tube
$]X[_{\P}$ of $X$ in $\P_K$ is the admissible open subset defined by
$$ ]X[_{\P}\ = \{ x \in \P_K\ |\ \forall i, |f_i(x)| < 1 \},$$
where, for a point $x \in \P_K$ corresponding to a maximal ideal 
$\mathfrak{m} \subset A_K$, $|f_i(x)|$ is the absolute value of the class 
of $f_i$ in the residue field $K(x) = A_K/\mathfrak{m}$. When $\P$ is 
not affine, one can choose a covering of $\P$ by open affine subsets 
$\U_i$, and the generic fibres $\U_{i, K}$ of the $\U_i$ provide an 
admissible covering of $\P_K$. Then the tube $]X[_{\P}$ can be defined 
by gluing the tubes $]X\cap\U_i[_{\U_i} \subset \U_{i,K}$.

\subsection{}\label{CohRig}
We introduce now the $\sO_{\P}$-algebra $\sA_{X,\P}$ of analytic
functions on $]X[_{\P}$, the de Rham cohomology of which defines rigid
cohomology.

The specialization morphism $\spm : \P_K\to\P$ maps $]X[_{\P}$ to $X$.
Thus, one can define a sheaf of $\sO_{\P,K}$-algebras supported on $X$
by setting
$$ \sA_{X,\P} = \spm_{\ast}\sO_{]X[_{\P}}.$$
The differentiation of analytic functions on $]X[_{\P}$ endows
$\sA_{X,\P}$ with a canonical integrable connection, allowing to
define the de Rham complex $\sA_{X,\P}\otimes\Omd{\P}$. For any
affine open subset $\U \subset \P$ and any $j$, we have
$\Gamma(\U_K,\Omega^j_{\P_K}) \cong \Gamma(\U, \Omega^j_{\P})\otimes K$.
It follows that there is a canonical isomorphism of complexes
$$ \sA_{X,\P}\otimes_{\sO_{\P}} \Omd{\P} \cong
\spm_{\ast}(\Omd{]X[_{\P}}).$$
Moreover, the above description of $]X[_{\P}$ in the affine case
shows that the inverse image of an affine open subset $\U \subset
\P$ is quasi-Stein in the sense of Kiehl \cite[Definition 2.3]{Ki}.
Therefore, Kiehl's vanishing theorem for coherent analytic sheaves
implies that $R^i\spm_{\ast}(\Omega^j_{]X[_{\P}}) = 0$ for all $j$ and
all $i \geq 1$, and we obtain
$$ \sA_{X,\P}\otimes_{\sO_{\P}} \Omd{\P} \cong
\R\spm_{\ast}(\Omd{]X[_{\P}})$$
in the derived category $\Db(X,K)$.

If $X \inj \P$ and $X \inj \P'$ are two closed immersions of $X$ into 
smooth formal $W$-schemes, there exists in $\Db(X,K)$ a canonical 
isomorphism \cite[1.5]{Be1}
\eq{IndImm}{ \sA_{X,\P}\otimes \Omd{\P} \cong \sA_{X,\P'}\otimes 
\Omd{\P'}. }
Thus, up to canonical isomorphism, the complex $\sA_{X,\P}\otimes
\Omd{\P}$ does not depend in $\Db(X,K)$ on the choice of the embedding
$X \inj \P$. In particular, the de Rham cohomology of $]X[_{\P}$ does
not depend on this choice. When $X$ is a proper $k$-scheme and can be
embedded in a smooth formal scheme $\P$ as above, its rigid cohomology
is defined by
\ga{defCohRig}{\R\Gamma\rig(X/K) = \R\Gamma(]X[_{\P}, \Omd{]X[_{\P}}) \cong 
\R\Gamma(X, \sA_{X,\P}\otimes\Omd{\P}).}
If such an embedding does not exist, one can choose a covering of $X$
by affine open subsets $X_i$, and, for each $i$, a closed immersion of
$X_i$ in a smooth affine formal scheme $\P_i$ over $W$. Using the
diagonal immersions for finite intersections $X_{i_0,\ldots,i_r} =
X_{i_0} \cap \ldots \cap X_{i_r}$, and the corresponding algebras
$\sA_{X_{i_0,\ldots,i_r},\P_{i_0}\times\cdots\times\P_{i_r}}$, one can
build a \v{C}ech-de Rham double complex. When $X$ is proper, the
cohomology of the associated total complex defines rigid cohomology.

\subsection{}\label{CohRigc}
The previous constructions can be extended as follows to define rigid 
cohomology with compact supports for separated $k$-schemes of finite 
type.

Let $U$ be such a scheme, $U \inj X$ an open immersion of $U$ in a
proper $k$-scheme $X$, $Y = X \setminus U$. Assume that there exists a
closed immersion $X \inj \P$ of $X$ into a smooth formal scheme $\P$,
and denote by $u : ]Y[_{\P} \inj ]X[_{\P}$ the inclusion of the tube
of $Y$ into the tube of $X$. Then $]Y[_{\P}$ is an admissible open
subset of $]X[_{\P}$, and, by construction \cite{BeL}, the rigid
cohomology with compact supports of $U$ is defined by
$$ \R\Gamma\rigc(U/K) = \R\Gamma(]X[_{\P}, (\Omd{]X[_{\P}} \to
u_{\ast}(\Omd{]Y[_{\P}}))\tot),$$
where the subscript t denotes the total complex associated to a double 
complex. Using the algebras $\sA_{X,\P}$ and $\sA_{Y,\P}$ defined in 
\ref{CohRig}, we can rephrase this definition as
\ga{defCohRigc}{\R\Gamma\rigc(U/K) = \R\Gamma(X, 
(\sA_{X,\P}\otimes_{\sO_{\P}}\Omd{\P} \to 
\sA_{Y,\P}\otimes_{\sO_{\P}}\Omd{\P})\tot).} 
This cohomology only depends on $U$, and has the usual properties of 
cohomology with compact supports.

If there is no embedding of a compactification $X$ of $U$ in a smooth 
formal $W$-scheme, the rigid cohomology with compact supports of $U$ can 
still be defined using coverings by affine open subsets and \v{C}ech 
complexes as in \ref{CohRig}.

\begin{prop}\label{AXPtoWitt}
Let $\P$ be a smooth and separated formal scheme over $W$, with
special fibre $P$, and let $X$ be a closed subscheme of $P$.

\romain The datum of a $\sigma$-semilinear lifting $F : \P \to \P$ of
the absolute Frobenius endomorphism of $P$ defines a ring homomorphism
\eq{DefAXPtoWitt}{ \sA_{X,\P} \to W\sO_{X,K}, }
functorial in $(X, \P, F)$. 

\romain Without assumption on the existence of $F$, there exists in
$\Db(X, K)$ a morphism
\eq{DefAXPdRtoWitt}{ a_{X,\P} : \sA_{X,\P}\otimes_{\sO_{\P}}\Omd{\P} \to
W\sO_{X,K}, }
functorial in $(X, \P)$, equal to the composed morphism 
$$ \sA_{X,\P}\otimes_{\sO_{\P}}\Omd{\P} \to \sA_{X,\P} 
\xrightarrow{\eqref{DefAXPtoWitt}} W\sO_{X,K} $$ 
whenever there exists a lifting of the Frobenius endomorphism on $\P$,
and compatible with the canonical isomorphism \eqref{IndImm}
$$ \sA_{X,\P}\otimes_{\sO_{\P}}\Omd{\P} \cong 
\sA_{X,\P'}\otimes_{\sO_{\P'}}\Omd{\P'} $$ 
for two closed immersions of $X$ into smooth and separated formal
schemes $\P$ and $\P'$.

\romain Let $U$ be an open subset of $X$, $Y = X \setminus U$, and
$\sI \subset \sO_X$ a coherent ideal such that $V(\sI) = Y$. Without
assumption on the existence of $F$, there exists in $\Db(X,K)$ a
morphism
\eq{DefAXYPdRtoWitt}{ a_{U,X,\P} : (\sA_{X,\P}\otimes_{\sO_{\P}}\Omd{\P} 
\to \sA_{Y,\P}\otimes_{\sO_{\P}}\Omd{\P})\tot \to 
W\sI_K, }
functorial with respect to morphisms $(U',X',\P') \to (U,X,\P)$ such 
that $U' \to U$ is proper, equal to \eqref{DefAXPdRtoWitt} when 
$U = X$, and compatible with the canonical isomorphisms 
$$ (\sA_{X,\P}\otimes_{\sO_{\P}}\Omd{\P} \to 
\sA_{Y,\P}\otimes_{\sO_{\P}}\Omd{\P})\tot \cong 
(\sA_{X,\P'}\otimes_{\sO_{\P'}}\Omd{\P'} \to 
\sA_{Y,\P'}\otimes_{\sO_{\P'}}\Omd{\P'})\tot $$ 
for two closed immersions of $X$ into smooth formal schemes $\P$ and $\P'$. 
\end{prop}

\begin{proof}
Let $\sJ \subset \sO_{\P}$ be the ideal defining $X$, $\sP(\sJ)$ the
divided power envelope of $\sJ$, with compatibility with the natural
divided powers of $p$, and $\sPh(\sJ)$ its $p$-adic completion. 
We recall first from \cite[1.9]{Be1} that there exists a functorial
homomorphism of $\sO_{\P}$-algebras
\eq{DefAXPtoPD}{ \sA_{X,\P} \to \sPh(\sJ)_K. }
Let $P_n$ be the reduction of $\P$ mod $p^n$, and let $\sJ_n =
\sJ\sO_{P_n}$ be the ideal of $X$ in $P_n$. Thanks to the
compatibility condition with the divided powers of $p$, there is a
canonical isomorphism
$$ \sPh(\sJ)/p^n\sPh(\sJ) \xrightarrow{\ \sim\ \,} \sP(\sJ_n), $$ 
where $\sP(\sJ_n)$ is the divided power envelope of $\sJ_n$ with the 
same compatibility condition, and we obtain 
$$ \sPh(\sJ) \xrightarrow{\ \sim\ \,} \varprojlim_n \sP(\sJ_n). $$ 
Thus, it suffices to define a compatible family of functorial ring
homomorphisms
\eq{PDtoWitt}{ \sP(\sJ_n) \to W_n\sO_X } 
to obtain a morphism \eqref{DefAXPtoWitt}. 

Let us assume that $\P$ is endowed with a lifting $F$ of the Frobenius
morphism. As $\sO_{\P}$ is $p$-torsion free, the homomorphism $F :
\sO_{\P} \to \sO_{\P}$ defines a section $s_F : \sO_{\P} \to
W\sO_{\P}$ of the reduction homomorphism $W\sO_{\P} \to \sO_{\P}$,
characterized by the fact that $w_i(s_F(x)) = F^i(x)$ for any $x \in
\sO_{\P}$ and any ghost component $w_i$ \cite[0 1.3]{I}. Composing
$s_F$ with the homomorphisms $W\sO_{\P} \to W\sO_P \to W_n\sO_P$ and
factorizing, we obtain for all $n \geq 1$ a homomorphism
$$ \sO_{P_n} \to W_n\sO_P \to W_n\sO_X, $$ 
which maps $\sJ_n$ to $VW_{n-1}\sO_X \subset W_n\sO_X$. The ideal
$VW_{n-1}\sO_X$ has a natural structure of divided power ideal
(compatible with the divided powers of $p$), defined by $(Vx)^{[i]} =
(p^{i-1}/i!)V(x^i)$ for all $i \geq 1$. Therefore, this homomorphism
factors through a homomorphism $\sP(\sJ_n) \to W_n\sO_X$. This 
provides the compatible family of homomorphisms defining 
\eqref{DefAXPtoWitt}, and it is clear that the homomorphism obtained 
in this way is functorial in $(X, \P, F)$. 

By composition with the augmentation morphism
$\sA_{X,\P}\otimes_{\sO_{\P}}\Omd{\P} \to \sA_{X,\P}$, we obtain a
morphism of complexes 
$$ a_{X,\P,F} : \sA_{X,\P}\otimes_{\sO_{\P}}\Omd{\P}
\to W\sO_{X,K}, $$
which is still functorial in $(X,\P,F)$. If $X \inj \P$
and $X \inj \P'$ are two closed immersions in smooth formal schemes
endowed with liftings of Frobenius $F$ and $F'$, we can endow $\P'' =
\P \times_W \P'$ with $F\times F'$, and embed $X$ diagonally into
$\P''$. Applying the functoriality of this construction to the two
projections, we obtain a commutative diagram
\eq{CompAPXtoWitt}{ \xymatrix{
\sA_{X,\P}\otimes_{\sO_{\P}}\Omd{\P}\ar[d]^{a_{X,\P,F}}\ar[r] & 
\sA_{X,\P''}\otimes_{\sO_{\P''}}\Omd{\P''}\ar[d]^{a_{X,\P'',F''}} &
\sA_{X,\P'}\otimes_{\sO_{\P'}}\Omd{\P'}\ar[l]\ar[d]^{a_{X,\P',F'}}\\
W\sO_{X,K} & W\sO_{X,K}\ar@{=}[l]\ar@{=}[r] & W\sO_{X,K}
} } 
in which the morphisms of the top row are quasi-isomorphisms. In
$\Db(X,K)$, the composed isomorphism
$\sA_{X,\P}\otimes_{\sO_{\P}}\Omd{\P} \cong
\sA_{X,\P'}\otimes_{\sO_{\P'}}\Omd{\P'}$ is the canonical isomorphism
\eqref{IndImm} \cite[1.5]{Be1}, which shows the compatibility asserted
in (ii). In the particular case where $\P' = \P$, this composed
isomorphism is the identity, and we obtain that, as a morphism of
$\Db(X,K)$, the morphism $a_{X,\P,F}$ does not depend on the choice of
$F$. A similar argument shows that, in $\Db(X,K)$, it is functorial
with respect to morphisms $(X',\P') \to (X,\P)$ without compatibility
with Frobenius liftings. Thus, if we define $a_{X,\P}$ to be the image
in $\Db(X,K)$ of $a_{X,\P,F}$, assertion (ii) is true when there
exists a lifting of Frobenius on $\P$.

In the general case, we can choose an affine covering $\P_i$ of $\P$, 
and a lifting of Frobenius $F_i$ on each $\P_i$. Then, if 
$X_{i_0,\ldots,i_n}=X\cap \P_{i_0} \cap \ldots \cap \P_{i_n}$ and 
$j_{i_0,\ldots,i_n}$ denotes its inclusion in $X$, the complex 
$\sA_{X,\P}\otimes_{\sO_{\P}}\Omd{\P}$ is quasi-isomorphic to the 
total complex associated to the double complex
\begin{multline}\label{CechAXP} \prod_i
j_{i\,\ast}(\sA_{X_i,\P_i}\otimes\Omd{\P_i}) \to \cdots \\
\to \prod_{i_0,\ldots,i_n} j_{i_0,\ldots,i_n\,\ast}
(\sA_{X_{i_0,\ldots,i_n},\P_{i_0}\times\cdots\times\P_{i_n}}
\otimes \Omd{\P_{i_0}\times\cdots\times\P_{i_n}}) \to \cdots, 
\end{multline}
while $W\sO_{X,K}$ is quasi-isomorphic to the \v{C}ech resolution
\eq{CechWOX}{ \prod_i j_{i\,\ast}(W\sO_{X_i,K}) \to \cdots \to 
\prod_{i_0,\ldots,i_n} j_{i_0,\ldots,i_n\,\ast}
(W\sO_{X_{i_0,\ldots,i_n},K}) \to \cdots. }
Then we can define $a_{X,\P}$ as the image in $\Db(X,K)$ of the 
morphism defined by the collection of all $a_{X_{i_0,\ldots,i_n},
\P_{i_0}\times\cdots\times\P_{i_n},F_{i_0}\times\cdots\times 
F_{i_n}}$, and assertion (ii) is verified as above.

Let $U$ be an open subset of $X$, and assume again that there exists a
lifting $F$ of the Frobenius endomorphism on $\P$. Then the previous
constructions can be applied both to $X$ and to $Y = X \setminus U$.
By functoriality, they provide a morphism of complexes
$$ (\sA_{X,\P}\otimes_{\sO_{\P}}\Omd{\P} 
\to \sA_{Y,\P}\otimes_{\sO_{\P}}\Omd{\P})\tot \to 
(W\sO_{X,K} \to W\sO_{Y,K}). $$ 
As the length $1$ complex $W\sO_{X,K} \to W\sO_{Y,K}$ is a resolution
of $W\sI_K$, we can define $a_{U,X,\P}$ as the image of
this morphism in $\Db(X,K)$. If $Y = \emptyset$, then we simply obtain
\eqref{DefAXPdRtoWitt}. To check the functoriality, we may assume
that $U'$ is dense in $X'$, because $W\sI'_K$ does not change if we
replace $X'$ by the closure of $U'$ in $X'$, by \ref{IndComp}, and
similarly $(\sA_{X,\P}\otimes_{\sO_{\P}}\Omd{\P} \to
\sA_{Y,\P}\otimes_{\sO_{\P}}\Omd{\P})\tot$ does not change up to
canonical isomorphism in $\Db(X,K)$, by the basic properties of rigid
cohomology. Then the properness of $U' \to U$ allows to assume that
$Y' = X' \setminus U'$ is defined by $\sI' = \sI\sO_{X'}$, and the
functoriality is clear. Finally, the compatibility with the canonical
isomorphism for two embeddings of $X$ into smooth formal schemes
results from the same assertion for the morphisms
\eqref{DefAXPdRtoWitt} relative to $X$ and $Y$.

When $F$ cannot be lifted to $\P$, one can proceed with \v{C}ech 
coverings as in (ii) to define \eqref{DefAXYPdRtoWitt}, and the same 
properties hold.
\end{proof}

\noindent\textit{Remark}. If $X$ is quasi-projective, one can always
find a closed immersion of $X$ into a smooth formal scheme $\P$
endowed with a lifting of Frobenius $F$, since it suffices to choose
for $\P$ an open subscheme of a projective space, endowed with the
endomorphism induced by some lifting of the Frobenius endomorphism of
the projective space.

\begin{thm}\label{ConstRigtoWitt}
Let $U$ be a separated $k$-scheme of finite type. In $\Db(K)$, there
exists a functorial morphism
\eq{DefRigtoWitt}{ a_U : \R\Gamma\rigc(U/K) \to
\R\Gamma\cpt(U,W\sO_{U,K}), }
equal, for a proper and smooth $k$-scheme $X$, to the canonical 
morphism
\eq{DefCristoWitt}{ \xymatrix{
\R\Gamma\rig(X/K) \ar[r]^-{\sim} & 
\R\Gamma\cris(X/W)_K \ar[r]^-{\sim} &  
\R\Gamma(X, \WOmd{}{X})_K \ar[d] \\
& & \R\Gamma(X, W\sO_X)_K. 
} }
\end{thm}

In \eqref{DefCristoWitt}, the first isomorphism is the comparison
isomorphism between rigid and crystalline cohomologies
\cite[1.9]{Be1}, the second one is the comparison isomorphism between
crystalline and de Rham-Witt cohomologies \cite[II, (1.3.2)]{I}, and
the third map is defined by the augmentation morphism of the de
Rham-Witt complex.

\begin{proof}
Let $U \inj X$ be an open immersion of $U$ in a proper $k$-scheme. Let
us assume first that there exists a closed immersion $X \inj \P$ of
$X$ into a smooth formal scheme $\P$. To define $a_U$,
we apply the functor $\R\Gamma(X, -)$ to the morphism 
$a_{U,X,\P}$ defined in \eqref{DefAXYPdRtoWitt}. Because of the 
compatibility property of \ref{AXPtoWitt} (iii), we obtain in this 
way a morphism of $\Db(K)$ which does not depend, up to canonical 
isomorphism, on the choice of $\P$. Using the functoriality of 
$a_{U,X,\P}$, it is easy to check that it does not depend either on 
the choice of the compactification $X$ of $U$, and that it is 
functorial with respect to $U$. 

When $U = X$ is a proper and smooth $k$-scheme, and $X$ is embeddable
in a smooth formal scheme $\P$ with a lifting $F$ of the Frobenius
endomorphism, the homomorphisms \eqref{PDtoWitt} used to construct
$a_{X,\P,F}$ are the homomorphisms used in \cite[II, (1.1.5)]{I} to
construct the morphism of complexes $\sPh(\sJ)\otimes\Omd{\P} \to
\WOmd{}{X}$ which defines the isomorphism between crystalline and de
Rham-Witt cohomologies. Therefore, the comparison isomorphisms
which appear in \eqref{DefCristoWitt} sit in a commutative diagram
$$ \xymatrix{
\R\Gamma(X, \sA_{X,\P}\otimes\Omd{\P}) \ar[d]\ar[r]^-\sim &
\R\Gamma(X, \sPh(\sJ)\otimes\Omd{\P,K}) \ar[d]\ar[r]^-\sim &
\R\Gamma(X, \WOmd{}{X,K})\ar[d] \\
\R\Gamma(X, \sA_{X,\P}) \ar[r] &
\R\Gamma(X, \sPh(\sJ)_K)\ar[r] &
\R\Gamma(X, W\sO_{X,K}),
} $$
in which the vertical maps are induced by the augmentation morphisms
of the complexes appearing in the upper row. The equality of the
morphisms \eqref{DefRigtoWitt} and \eqref{DefCristoWitt} follows
immediately.

In the general case where there is no smooth embedding $\P$ of $X$, or
no Frobenius lifting on $\P$, one can again use a covering of $X$ by
affine open subsets $X_i$, and closed immersions $X_i \inj \P_i$ in
smooth formal schemes $\P_i$ endowed with liftings of Frobenius, to
construct \v{C}ech complexes as in \eqref{CechAXP} and
\eqref{CechWOX}. One can then use on each intersection of the covering
the corresponding morphism \eqref{DefAXYPdRtoWitt} and define in this
way a morphism between the two \v{C}ech complexes. Applying
$\R\Gamma(X, -)$, one gets in $\Db(K)$ the morphism $a_U$. In the
proper and smooth case, the fact that each map in the above diagram
can be defined through \v{C}ech complexes shows that $a_U$ is still
equal to \eqref{DefCristoWitt}.
\end{proof}

The following verifications are left as an exercise for the reader.

\begin{prop}\label{CompataX}
The morphism $a_U$ defined by Theorem \ref{ConstRigtoWitt} satisfies 
the following compatibility properties.

\romain The diagram
\eq{CompRed}{ \xymatrix{
\R\Gamma\rigc(U/K) \ar[r]^-\sim\ar[d]^-{a_U} & 
\R\Gamma\rigc(U_\red/K) \ar[d]^-{a_{U_\red}} \\
\R\Gamma\cpt(U,W\sO_{U,K}) \ar[r]^-\sim &
\R\Gamma\cpt(U,W\sO_{U_\red,K})
} }
is commutative.

\romain If $V \subset U$ is an open subset, and $T = U \setminus V$, 
the morphisms $a_V$, $a_U$ and $a_T$ define a morphism of exact 
triangles  
\ml{CompExSeq}{ \xymatrix{
\R\Gamma\rigc(V/K) \ar[r]\ar[d]^-{a_V} & 
\R\Gamma\rigc(U/K) \ar[r]\ar[d]^-{a_U} & 
\R\Gamma\rigc(T/K) \ar[r]^-{+1}\ar[d]^-{a_T} & \ \, \\
R\Gamma\cpt(V,W\sO_{V,K}) \ar[r] &
R\Gamma\cpt(U,W\sO_{U,K}) \ar[r] &
R\Gamma\cpt(T,W\sO_{T,K}) \ar[r]^-{+1} & \raisebox{-2mm}{\ .}
} }
\end{prop}

\bigskip
\section{Proof of the main theorem}
\label{ProofMainTh}
\medskip

We prove here that the morphism $a_U$ constructed in Theorem
\ref{ConstRigtoWitt} yields an identification of the slope $< 1$
subspace of rigid cohomology with compact supports with Witt vector
cohomology with compact supports, thus completing the proof of Theorem
\ref{SlopeHc}. 

The next lemma will allow us use descent techniques to study the
morphism $a_U$. We follow here the method of Chiarellotto and Tsuzuki
to construct embeddings of a simplicial scheme in a smooth simplicial
formal scheme (\cf \cite[11.2]{CT}, \cite[7.3]{Ts2}).

\begin{lem}\label{EmbedResol}
Let $X$ be a separated reduced $k$-scheme of finite type, and let $X
\inj \P$ be a closed immersion into a smooth formal $W$-scheme endowed
with a Frobenius lifting $F$. Let $N \in \N$ be a fixed integer. There
exists a proper hypercovering $X\lbul$ of $X$, a $\P$-augmented
simplicial formal $W$-scheme $\P\lbul$, endowed with a
$\sigma$-semilinear endomorphism $F\lbul$ lifting the absolute
Frobenius endomorphism of its special fibre, and a morphism of
simplicial schemes $X\lbul \to \P\lbul$ above $X \inj \P$, such that
the following conditions are satisfied:

\romain For all $n$, the morphism $X_n \to \P_n$ is a closed 
immersion, and the projection $\P_n \to \P$ commutes with the 
Frobenius liftings.

\romain The $N$-truncated simplicial scheme $\sk^X_N(X\lbul)$ is a 
$N$-truncated simplicial resolution of $X$.

\romain The canonical morphisms $X\lbul \to \cosk^X_N(\sk^X_N(X\lbul))$ 
and $\P\lbul \to \cosk^{\P}_N(\sk^{\P}_N(\P\lbul))$ are isomorphisms. 

\romain For all $n$, $\P_n$ is smooth over $\P$, and the canonical
morphism $\P_n \to \cosk^{\P}_{n-1}(\sk^{\P}_{n-1}(\P\lbul))_n$ is
smooth.
\end{lem}

\begin{proof} 
Using de Jong's theorem, we can find an $N$-truncated simplicial
resolution of $X$ (\cf \ref{resolutions}). We define $X\lbul$ as the
$N$-coskeleton over $X$ of this $N$-truncated resolution. Thus
$X\lbul$ is a proper hypercovering of $X$, $\sk^X_N(X\lbul)$ is our
initial $N$-truncated simplicial resolution of $X$, and the morphism
$X\lbul \to \cosk^X_N(\sk^X_N(X\lbul))$ is an isomorphism.

For $n \in \N$, let $[n]$ denote the ordered set $\{0,\ldots,n\}$. Let
$\Delta$ be the category which has the sets $[n]$ as objects, and the
set of non decreasing maps $[m] \to [n]$ as set of morphisms from
$[m]$ to $[n]$. Since the $X_m$ are quasi projective over $k$, we can
choose for each $m \leq N$ a closed immersion $i_m : X_m \inj \W_m$,
where $\W_m$ is a smooth formal $W$-scheme endowed with a lifting of
Frobenius $F_m$. From $\W_m$, we construct a simplicial complex of smooth
formal $W$-schemes $\Gamma_n(\W_m)$ by setting for all $n \geq 0$
\eq{DefGamma}{ \Gamma_n(\W_m) = \prod_{\gamma:[m]\to[n]}\W_{m,\gamma}, }
where the product is taken over all morphisms $\gamma:[m]\to[n]$ in
$\Delta$, and $\W_{m,\gamma} = \W_m$ for all $\gamma$. If $\eta :
[n']\to[n]$ is a morphism in $\Delta$, the corresponding morphism
$\Gamma_n(\W_m) \to \Gamma_{n'}(\W_m)$ is such that, for any $\gamma'
: [m] \to [n']$ in $\Delta$, its composition with the projection
$\Gamma_{n'}(\W_m) \to \W_{m,\gamma'} = \W_m$ is the projection of
$\Gamma_n(\W_m)$ on $\W_{m,\eta\circ\gamma'} = \W_m$. We can then
define a Frobenius lifting on the simplicial formal scheme
$\Gamma\lbul(\W_m)$ as being the product morphism defined by $F_m$ on
each $\Gamma_n(\W_m)$; we still denote it by $F_m$.

The immersion $i_m$ defines a morphism of simplicial schemes $X\lbul
\to \Gamma\lbul(\W_m)$ as follows: for each $n \geq 0$, the
composition of $X_n \to \Gamma_n(\W_m)$ with the projection of index
$\gamma$ is the morphism $X_n \to X_m \inj \W_m$, where the first
morphism is the morphism of $X\lbul$ defined by $\gamma$. We can now
define a simplicial formal $W$-scheme $\P\lbul$ augmented towards $\P$
by setting
\eq{DefPP}{ \P\lbul = \P \times_{\Spf(W)} \prod_{0 \leq m \leq N} 
\Gamma\lbul(\W_m). }
We define a morphism of simplicial schemes $X\lbul \to \P\lbul$ using
the composed morphism $X\lbul \to X \to \P$ and the family of
morphisms $X\lbul \to \Gamma\lbul(\W_m)$ defined above. By
construction, this morphism is compatible with the immersion $X \inj
\P$ via the augmentation morphisms. If we endow $\P\lbul$ with the
Frobenius lifting $F \times \prod_m F_m$, then the augmentation
$\P\lbul \to \P$ commutes with the Frobenius liftings.

To check the remaining properties, we observe first that $\W\lbul
\xrightarrow{\ \sim\ \,} \cosk^{\P}_N(\sk^{\P}_N(\W\lbul))$ thanks to
\cite[11.2.5]{CT}. For each $m \leq N$, the morphism $X_m \to
\Gamma_m(\W_m)$ is a closed immersion since the chosen immersion $X_m
\inj \W_m$ is one of the factors. Therefore $X_n \to \P_n$ is a closed
immersion for all $n \leq N$. Then, thanks to the previous property,
it follows from the construction of coskeletons that $X_n \to \P_n$ is
a closed immersion for all $n$. Finally, each $\P_n$ is smooth over
$\P$ by construction, and the fact that $\P_n$ is smooth over
$(\cosk^{\P}_{n-1}(\sk^{\P}_{n-1}(\P\lbul))_n$ for all $n$ follows
from \cite[7.3.3]{Ts2}.
\end{proof}

\subsection{}\label{DemSlopeHc}\textit{Proof of Theorem \ref{SlopeHc}.}
We first recall the definition of the subspaces with prescribed
slopes. Let $\kbar$ be an algebraic closure of $k$, and let $\Kbar_0 =
\Frac(W(\kbar))$. If $E$ is a finite dimensional $K$-vector space
endowed with a $\sigma$-semilinear automorphism $F$, let $\Ebar =
\Kbar_0\otimes_K E$, and $\Fbar = \sigma_{\Kbar_0}\otimes F$. The
Dieudonn\'e-Manin theorem provides for each $\lambda \in \Q$ a maximal
$\Fbar$-stable subspace $\Ebar^{\lambda}$ on which $\Fbar$ has purely
slope $\lambda$, and asserts that $\bigoplus_{\lambda}\Ebar^{\lambda}
\cong E$. As these subspaces are invariant under the action of
$\Gal(\kbar/k)$, there exists a unique $F$-stable subspace
$E^{\lambda} \subset E$ such that $\Kbar_0 \otimes_K E^{\lambda}
\xrightarrow{\ \sim\ \,} \Ebar^{\lambda}$, and
$\bigoplus_{\lambda}E^{\lambda} \cong E$. The subspace $E^{\lambda}$
is the \textit{slope $\lambda$ subspace} of $E$, and, for any $\rho
\in \R$, we define the \textit{slope $< \rho$ subspace} of $E$ by
$$ E^{< \rho} = \bigoplus_{\lambda < \rho} E^{\lambda}. $$ 
These spaces are functorial with respect to $(E, F)$, and define 
exact functors.

Let $X$ be a separated $k$-scheme of finite type. The morphism 
\eqref{DefRigtoWitt} induces canonical homomorphisms
$$ a^q_X : H^q\rigc(X/K) \to H^q\cpt(X, W\sO_{X,K}) $$ 
between the cohomology spaces. By functoriality, $a_X^q$ commutes to
the Frobenius actions on both sides. The spaces $H^q\rigc(X/K)$ are
finite dimensional \cite[3.9~(i)]{Be1}, and their Frobenius
endomorphism is an automorphism \cite[2.1]{ELS}. Therefore they have a
slope decomposition, and $a^q_X$ induces a homomorphism
$$ b^q_X : H^q\rigc(X/K)^{<1} \to H^q\cpt(X, W\sO_{X,K}) $$ 
which is the homomorphism \eqref{sl3} of Theorem \ref{SlopeHc}.

To prove that $b^q_X$ is an isomorphism, we first observe that this is
true when $X$ is proper and smooth over $k$. Indeed, $a^q_X$ is then
induced by the morphism \eqref{DefCristoWitt}, and the isomorphisms
which enter in the definition of \eqref{DefCristoWitt} are compatible
with Frobenius, hence induce isomorphisms on the slope $< 1$ subspaces 
of the cohomology spaces. The theorem is then a consequence of the slope
decomposition for de Rham-Witt cohomology \cite[II,~(3.5.2)]{I}.

In the general case, we prove Theorem \ref{SlopeHc} by induction on
$\dim(X)$. The diagram \eqref{CompRed} allows us to assume that $X$ is
reduced. In particular, the theorem holds for $\dim(X) = 0$ thanks to
the previous remark. Let us assume that all $b^q_X$ are isomorphisms
when $\dim(X) < d$. If $U \subset X$ is an open subset of dimension
$d$ of a proper and smooth $k$-scheme $X$, its closure $\Ubar$ is
proper and smooth over $k$, and $\dim(\Ubar\setminus U) < d$. As all
$b^q_{\Ubar}$ are isomorphisms, the induction hypothesis and the
commutativity of \eqref{CompExSeq} imply that all $b^q_U$ are
isomorphisms. Thus the theorem also holds for any open subset of
dimension $d$ of a proper and smooth $k$-scheme.

Let $X$ be an arbitrary reduced and separated $k$-scheme of dimension
$d$. Using again \eqref{CompExSeq} and the induction hypothesis, we
may replace $X$ by an arbitrarily small neighbourhood $U$ of the
generic points of its irreducible components. As the theorem holds for
a scheme if and only if it holds for each of its connected components,
we may shrink even more and assume that $U$ is irreducible and affine.
Let now $U \inj X$ be an open immersion of $U$ in a projective
$k$-scheme $X$, such that $U$ is dense in $X$. Using again the same
argument, it suffices to prove the theorem for $X$.

Since $X$ is projective, we can choose a closed immersion of $X$ into 
a smooth formal $W$-scheme $\P$ endowed with a Frobenius lifting $F$. 
Let $N$ be an integer such that $N \geq 2d$. Then we can choose a 
proper hypercovering $X\lbul$ of $X$, a $\P$-augmented simplicial 
formal scheme $\P\lbul$ endowed with a Frobenius lifting $F\lbul$, 
and a closed immersion of augmented simplicial schemes $X\lbul \inj 
\P\lbul$ above $X \inj \P$, so that the conditions of Lemma 
\ref{EmbedResol} are satisfied.

In particular, there is a dense open subset $U \subset X$ such that 
the restriction $U\lbul$ of $X\lbul$ above $U$ is an \'etale 
hypercovering of $U$, and it suffices to prove the theorem for $U$. 
Let $Y = X \setminus U$. Since each $\P_n$ is smooth over $W$, the 
complexes 
$(\sA_{X_n,\P_n}\otimes\Omd{\P_n} \to 
\sA_{Y_n,\P_n}\otimes\Omd{\P_n})\tot$ are defined. By functoriality, 
they define a complex of sheaves on $X\lbul$ \cite[(5.1.6)]{De1}, 
and we set
$$ H^{\ast}\rigc(U\lbul/K) := H^{\ast}(X\lbul, 
(\sA_{X\lbul,\P\lbul}\otimes\Omd{\P\lbul} \to 
\sA_{Y\lbul,\P\lbul}\otimes\Omd{\P\lbul})\tot). $$
As usual, there is a biregular spectral sequence relating the 
cohomology of the individual complexes on the $X_n$'s to the 
cohomology of the global complex on $X\lbul$, which can be written here
$$ E_1^{i,j} = H^j\rigc(U_i/K) \Rightarrow H^{i+j}\rigc(U\lbul/K). $$ 
By functoriality, this spectral sequence is endowed with a Frobenius 
automorphism. Moreover, the $E_1^{i,j}$ terms are finite dimensional 
$K$-vector spaces, so the spaces $H^{i+j}\rigc(U\lbul/K)$ are finite 
dimensional too. Since the subspace of slope $< 1$ is an exact 
functor on the category of finite dimensional $K$-vector spaces 
endowed with a $\sigma$-semilinear automorphism, we obtain a 
biregular spectral 
sequence
\eq{SpecSeqRig}{ E_1^{i,j} = H^j\rigc(U_i/K)^{<1} \Rightarrow 
H^{i+j}\rigc(U\lbul/K)^{<1}. }

Let $\sI \subset \sO_X$ be a coherent ideal such that $V(\sI) = Y$, 
and set $\sI_n = \sI\sO_{X_n}$, $Y_n = V(\sI_n) \subset X_n$. The 
sheaves $W\sO_{X_n}$, $W\sO_{Y_n}$ and $W\sI_n$ define sheaves 
$W\sO_{X\lbul}$, $W\sO_{Y\lbul}$ and $W\sI\lbul$ on $X\lbul$, with an 
exact sequence
$$ 0 \to W\sI\lbul \to W\sO_{X\lbul} \to W\sO_{Y\lbul} \to 0. $$ 
We define
$$ H^{\ast}\cpt(U\lbul, W\sO_{U\lbul\!,K}) = H^{\ast}(X\lbul, 
W\sI\lbul{}_{\!,K}), $$ 
and we obtain a biregular spectral sequence
\eq{SpecSeqWitt}{ E_1^{i,j} = H^j\cpt(U_i, W\sO_{U_i,K}) \Rightarrow 
H^{i+j}\cpt(U\lbul, W\sO_{U\lbul\!,K}). }

Because $F$ is an endomorphism of the simplicial scheme $\P\lbul$, the
homomorphisms \eqref{DefAXPtoWitt} define morphisms of sheaves of
rings $\sA_{X\lbul\!,\P\lbul} \to W\sO_{X\lbul\!,K}$ and
$\sA_{Y\lbul\!,\P\lbul} \to W\sO_{Y\lbul\!,K}$ on $X\lbul$, from which we
derive a morphism of complexes
$$ (\sA_{X\lbul\!,\P\lbul}\otimes\Omd{\P\lbul} \to 
\sA_{Y\lbul\!,\P\lbul}\otimes\Omd{\P\lbul})\tot \to 
(W\sO_{X\lbul\!,K} \to W\sO_{Y\lbul\!,K}). $$ 
Taking cohomology, we obtain homomorphisms
$$ b^q_{U\lbul} : H^q\rigc(U\lbul/K)^{<1} \to 
H^q\cpt(U\lbul, W\sO_{U\lbul\!,K}), $$ 
and a morphism of spectral sequences from \eqref{SpecSeqRig} to 
\eqref{SpecSeqWitt}. On the $E_1^{i,j}$ terms, this morphism is given 
by the homomorphisms $b^j_{U_i}$. 

Since the truncation $\sk^X_N(X\lbul)$ is a $N$-truncated resolution
of $X$, $X_n$ is a projective and smooth $k$-scheme for all $n \leq N$.
Then the theorem holds for the open subset $U_n \subset X_n$, which
has dimension $d$ since it is \'etale over $U$. Therefore, the
homomorphisms $b^j_{U_i}$ between the $E_1^{i,j}$ terms of
\eqref{SpecSeqRig} and \eqref{SpecSeqWitt} are isomorphisms for $i
\leq N$. It follows that $b^q_{U\lbul}$ is an isomorphism for $q \leq
N$.

As the augmentation $\P\lbul \to \P$ is compatible with the Frobenius 
liftings, we obtain for all $q$ a commutative square
\eq{CompbX}{ \xymatrix{
H^q\rigc(U\lbul/K)^{<1} \ar[r]^-{b^q_{U\lbul}} & 
H^q\cpt(U\lbul, W\sO_{U\lbul\!,K}) \\
H^q\rigc(U/K)^{<1} \ar[u] \ar[r]^-{b^q_U} & 
H^q\cpt(U, W\sO_{U,K}). \ar[u]
} }
The augmented simplicial schemes $X\lbul \to X$ and $Y\lbul \to Y$ are
proper hypercoverings, and the condition (iv) of \ref{EmbedResol}
implies that $(X\lbul,X\lbul,\P\lbul)$ and $(Y\lbul, Y\lbul,\P\lbul)$
are Pr-Pr-Sm-hypercoverings of $(X,X,\P)$ and $(Y,Y,\P)$ in the sense
of \cite[5.1.5]{Ts2}. Therefore, it follows from Tsuzuki's proper
descent theorem (\cite[2.1.3]{Ts1}, \cite[5.3.1]{Ts2}) that the left
vertical arrow is an isomorphism. On the other hand, the right
vertical arrow is an isomorphism thanks to Theorem \ref{descent}.
Since the source and target of $b^q_U$ are $0$ for $q > 2d$, and
$b^q_{U\lbul}$ is an isomorphism for $q \leq 2d$, this completes the
proof of Theorem \ref{SlopeHc}. \hfill $\Box$

\bigskip
\section{Applications and examples}
\label{Applications}
\medskip

In this section, we assume that $k$ is a finite field with $q = p^a$ 
elements, and we give applications of Theorem \ref{SlopeHc} to 
congruences mod $q$ for the number of rational points of some 
algebraic varieties. 

If $X$ is a separated $k$-scheme of finite type, we set for all $i$
$$ P_i(X, t) = \det(1 - t\phi|H^i\rigc(X/K)), $$ 
where $\phi = F^a$ is the $k$-linear Frobenius endomorphism of $X$. We
recall that, thanks to the Lefschetz trace formula \cite{ELS}, the
zeta function of $X$ can be expressed as
\eq{zeta}{ \zeta(X,t) = \prod_i P_i(X,t)^{(-1)^{i+1}}. } 

\subsection{}\label{Pless1}
We normalize the valuation on $K$ by setting $v_q(q) = 1$. If $P(t) =
\sum_i a_i t^i \in K[t]$, its Newton polygon is the graph of the
greatest convex function $s$ on $[0,\deg(P)]$ such that $s(i) \leq
v_q(a_i)$ for all $i$. If $P$ is irreducible, then its Newton polygon 
is a segment. For any polynomial $P(t) \in 1+tK[t]$, and any
$\lambda \in \Q$, we denote by $P^{\lambda}(t)$ the product of the
irreducible factors of $P$ whose Newton polygon is a segment of slope
$\lambda$, normalised by $P^{\lambda}(0)=1$. For any $\rho \in \R$, we 
set $P^{<\rho}(t) = \prod_{\lambda < \rho} P^{\lambda}(t)$. 

This definition extends by multiplicativity to rational fractions
$R(t) \in K(t)$ such that $R(0)=1$. In particular, it can be applied
to $\zeta(X,t) \in \Q(t) \subset K(t)$, for any separated $k$-scheme
of finite type $X$, and this defines the \textit{slope $< \rho$ factor}
$\zeta^{<\rho}(X,t)$ of $\zeta(X,t))$. By Manin's theorem relating over
finite fields the slopes with the eigenvalues of $\phi$, we obtain
$$ \det(1 - t\phi|H^i\rigc(X/K)^{<\rho}) = \det(1 - 
t\phi|H^i\rigc(X/K))^{<\rho}. $$ 

On the other hand, we define
\ga{}{ P^W_i(X,t) = \det(1 - t\phi|H^i\cpt(X, W\sO_{X,K})), \notag\\
\zeta^W(X,t) = \prod_i P^W_i(X,t)^{(-1)^{i+1}}. \notag} 
Then Theorem \ref{SlopeHc} implies Corollary \ref{zetaW}:
\ga{fPiW}{ \forall i,\ \ P^{<1}_i(X,t) = P^W_i(X,t), \\
\zeta^{<1}(X,t) = \zeta^W(X,t). \notag} 

We observe that these polynomials actually have coefficients in 
$\Z_p$:

\begin{prop}\label{CoefZp}
For any separated $k$-scheme of finite type $X$, and any $\lambda \in
\Q$, the polynomials $P_i^{\lambda}(X,t)$ and $P^W_i(X,t)$ belong to
$\Z_p[t]$.
\end{prop}

\begin{proof}
By construction, the polynomials $P_i(X,t)$ belong to $K[t]$. Let 
$\phi' : X' \to X'$ be the pull-back of $\phi$ by the absolute Frobenius 
endomorphism of $\Spec(k)$. By base change, we obtain
$$ \det(1 - t\phi'|H^i\rigc(X'/K)) = 
\sigma^{\ast}(P_i(X, T)). $$ 
On the other hand, the relative Frobenius $F_{X/k} : X \to X'$ 
commutes with $\phi$ and $\phi'$. As it induces an isomorphism 
$H^i\rigc(X'/K) \xrightarrow{\ \sim\ \,} H^i\rigc(X/K)$, it follows 
that 
$$ \det(1 - t\phi'|H^i\rigc(X'/K)) = P_i(X, T). $$ 
Thus $P_i(X, t)$ is invariant under $\sigma$, hence belongs to
$\Q_p[t]$. Moreover, since the Newton polygon of a polynomial does not
change by enlarging the field containing its coefficients, the
decomposition $P_i(X,t) = \prod_{\lambda} P_i^{\lambda}(X,t)$ is
defined in $\Q_p[t]$.

Finally, the slopes of the Frobenius action on $H^i\rigc(X/K)$ are
non-negative \cite[3.1.2]{CLS}, and this is equivalent by Manin's theorem
to the fact that the inverses of the roots of $P_i(X,t)$ are $p$-adic
integers. Therefore, each $P_i^{\lambda}(X,t)$ belongs to $\Z_p[t]$.

By \eqref{fPiW}, the statement for $P^W_i(X,t)$ follows (it can also 
be proved directly by the same argument).
\end{proof}

We recall the following well-known result \cite{Ax}, which links
congruences modulo some power of $q$ to the triviality of the
corresponding slope factor:

\begin{prop}\label{Divis}
Let $(N_r)_{r \geq 1}$ be a sequence of integers such that
$$\Phi(t) = 
\exp(\sum_{r \geq 1} N_r\frac{t^r}{r}) \ \in 1 + \Z[[t]],$$
and assume that $\Phi(t)$ is a rational function in $\Q(t)$. Then $\Phi(t)$ can
be written
$$ \Phi(t) = \frac{\prod_i(1-\alpha_i t)}{\prod_j(1-\beta_j t)}, $$
where $\alpha_i, \beta_j \in \overline{\Z}$ are such that $\alpha_i
\neq \beta_j$ for all $(i,j)$, and, for any integer $\kappa \geq 1$,
the following conditions are equivalent:

\romain For all $r \geq 1$, $N_r \equiv 0 \mod q^{\kappa r}$.

\romain For all $i$ and all $j$, $\alpha_i$ and $\beta_j$ are
divisible by $q^{\kappa}$ in $\overline{\Z}$.

\romain For all embeddings $\iota$ of $\overline{\Q}$ in
$\overline{\Q}_p$, all $i$ and all $j$, $v_q(\iota(\alpha_i)) \geq \kappa$,
$v_q(\iota(\beta_j)) \geq \kappa$.

\romain $\Phi^{<\kappa}(t) = 1$.
\end{prop}

As explained in the introduction, our first application of Theorem
\ref{SlopeHc} over finite fields will be a proof of a consequence of
Serre's conjecture on theta divisors.

\begin{prop}\label{GeneralSerre}
Let $X$ be a smooth projective variety of pure dimension $n$ over a
finite field. Let $D\subset X$ be an ample divisor, of complement
$U=X\setminus D$, and assume that the following conditions hold:

\alphab The injection $H^0(X, \omega_X)\to H^0(X, \omega_X(D))$ is an
isomorphism;

\alphab For all $i \geq 1$, the canonical homomorphism
$$ H^n(X, \sO_X) \xrightarrow{V^i} H^n(X, W_{i+1}\sO_X) $$
is injective.

Then $H^n\rigc(U/K)^{<1}\to H^n\rig(X/K)^{<1}$ is an isomorphism, and
\eq{zetaD}{ \zeta^{<1}(X, t) =
\zeta^{<1}(D, t) \cdot P^{<1}_n(X,t)^{(-1)^{n+1}}. }
\end{prop}

\begin{proof}
Let $\sI \cong \sO_X(-D)$ be the ideal of $D$ in $X$. Condition~a) is
\hfill\linebreak
equivalent by Serre duality to the condition
$$ H^n(X, \sI) \xrightarrow{\ \sim\ \,} H^n(X, \sO_X). $$
Thanks to condition~b), the rows of the commutative diagram
$$ \xymatrix@1@M=0.1cm@C=0.4cm@R=0.45cm{
0 \ar[r] & H^n(X, \sI) \ar[r]\ar[d]^{\wr} &
H^n(X, W_{i+1}\sI) \ar[r]\ar[d] & H^n(X, W_i\sI) \ar[r]\ar[d] & 0 \\
0 \ar[r] & H^n(X, \sO_X) \ar[r] &
H^n(X, W_{i+1}\sO_X) \ar[r] & H^n(X, W_i\sO_X) \ar[r] & 0
} $$
are exact, hence the homomorphism $H^n(X, W_i\sI) \to
H^n(X, W_i\sO_X)$ is an isomorphism for all $i$. Taking inverse limits
and tensoring with $K$, we obtain the isomorphism
\eq{TopCohU}{ H^n\cpt(U, W\sO_{U,K}) \xrightarrow{\ \sim\ \,} H^n(X,
W\sO_{X,K}). }
Thus Theorem \ref{SlopeHc} implies that $H^n\rigc(U/K)^{<1} \to
H^n\rig(X/K)^{<1}$ is also an isomorphism.

The multiplicativity of the zeta function shows that it suffices to
prove that
$$ \zeta^{<1}(U,t) = P^{<1}_n(X,t)^{(-1)^{n+1}}. $$
By \eqref{fPiW}, we have
$$ P^{<1}_i(U,t) = P^W_i(U,t) $$
for all $i$. By \eqref{TopCohU}, we have $P^W_n(U,t) = P^W_n(X,t) =
P^{<1}_n(X,t)$. On the other hand, $U$ is affine and smooth of pure
dimension $n$, hence $H^i\cpt(U, W\sO_{U,K})= 0$ for $i \neq n$ by
Corollary \ref{Vanish}. The proposition follows.
\end{proof}

\subsection{}\label{ProofSerre}\textit{Proof of Theorem \ref{ConjSerre}.}
Applying Proposition \ref{Divis}, it suffices to prove that, if
$\Theta$ and $\Theta'$ are two theta divisors in an abelian variety
$A$, then $\zeta^{<1}(\Theta,t) = \zeta^{<1}(\Theta',t)$. This will
follow from \eqref{zetaD} if we check that conditions a) and b) of
Proposition \ref{GeneralSerre} hold for a theta divisor in an abelian
variety. Since $A$ is an abelian variety, the dualizing bundle
$\omega_A$ is trivial, and, for a theta divisor $\Theta$, the
homomorphism $H^0(A, \sO_A) \to H^0(A, \sO_A(\Theta))$ is an
isomorphism \cite[III 16]{Mu}. Thus condition a) is satisfied. As for
condition b), the injectivity of the homomorphisms $H^n(A,\sO_A) \to
H^n(A, W_{i+1}\sO_A)$ follows from Serre's theorem on the vanishing of
Bockstein operations for abelian varieties in characteristic~$p$
(\cite[\S 1, 3]{Se} and \cite[Th\'eor\`eme 2]{Se2}). \hfill$\Box$

\begin{rmks}\label{RemonSerre} \hspace{-8mm} 
\romain In fact, Serre's conjecture is phrased more motivically: the
difference of the motives of $\Theta$ and $\Theta'$ should be
divisible in a suitable sense by the Lefschetz motive. What we show in
Theorem \ref{ConjSerre} is a finite field implication of the motivic
assertion. The motivic statement is clearly stronger than this
implication. On the other hand, it is hard to approach directly as it
deals with non-effective motives.

\romain The property of the Theorem is very special for theta
divisors. It is of course not true in general that two effective
divisors $D, D'$ with $h^0(\sO(D))=h^0(\sO(D'))$ carry the same number
of points modulo $q$: take for example $D=2\{0\}, D'=\{0\}+ \{
\infty\}$ on $\P^1$.

\romain Let us also remark that over the complex numbers, when $\Theta$ is
irreducible, one knows that its singularities are rational
\cite[Theorem 3.3]{EL}. A precise analogy of this over a finite field
isn't quite clear, as the notion of rational singularities itself
requires resolution of singularities. But it should be related to the
assertion of Corollary \ref{zetaW}.
\end{rmks}

According to Grothendieck-Deligne's philosophy of motives, Hodge type
over the field of complex numbers behaves the same way as
congruences for the number of rational points. In view of the previous
remarks, it is worth pointing out that Theorem \ref{ConjSerre} has the
following Hodge theoretic analogue (where $H^{\ast}(X)$ denotes the
classical complex cohomology of an algebraic scheme $X$ over $\C$):

\begin{prop} \label{PropHodge}
Let $A$ be an abelian variety of dimension $n$ over $\C$, and $\Theta
\subset A$ a theta divisor. Let $F$ be the Deligne's Hodge filtration
on $H^i(A)$ and $H^i(\Theta)$. Then the restriction map $\gr^0_F
H^i(A)\to \gr^0_F H^i(\Theta)$ is an isomorphism for $i \neq n$. In
particular, if $\Theta'$ is another theta divisor, $\dim \gr^0_F
H^i(\Theta)= \dim \gr^0_F H^i(\Theta')$ for all $i$.
\end{prop}

\begin{proof}
This is a simple example of application of \cite[Proposition 1.2]{E}.
Since $\Theta$ is ample, we know by the weak Lefschetz theorem that
the restriction map $H^i(A)\to H^i(\Theta)$ is an isomorphism for $i
\leq n-2$, and is injective for $i = n-1$. Also $\gr^0_F H^i(X) = 0$
for $i > \dim(X)$ for any separated $\C$-scheme of finite type $X$
\cite[Th. 8.2.4]{De1}. Thus all we have to prove is
\ga{IsomGr0}{\gr^0_F H^n_c(U) \xrightarrow{\sim}  \gr^0_F H^n(A),}  
with $U = A \setminus \Theta$. Let $\sigma: A'\to A$ be a birational
morphism such that $\sigma|_U$ is an isomorphism, $A'$ is smooth, and
$D = \sigma^{-1}(\Theta)_{\red}$ is a normal crossings divisor. Then one
has $\gr^0_F H_c^n(U) = H^n(A', \sO_{A'}(-D))$, and $\sigma^*: \gr^0_F
H^n(A) = H^n(A, \sO_A) \to \gr^0_F H^n(A') = H^n(A', \sO_{A'})$ is an
isomorphism. One the other hand, the composed morphism 
\ga{4.3}{H^n(A, \sO_A(-\Theta)) \xrightarrow{\sigma^*} H^n(A', 
\sO_{A'}(-D))
\xrightarrow{\iota} H^n(A', \sO_{A'}) \xleftarrow{\sim} H^n(A, \sO_A)}
is an isomorphism since $\Theta$ is a theta divisor, while $\sigma^*,
\iota$ are surjective for dimension reasons. Thus all maps in
\eqref{4.3} are isomorphisms, which in particular proves the
proposition.
\end{proof}

We now consider the case of intersections of hypersurfaces of small 
degrees.

\begin{prop}\label{CohAx}
Let $k$ be a perfect field of characteristic $p$, let $K = \Frac(W(k))$, 
and let $D_1, \ldots, D_r \subset \P^n_k$ be hypersurfaces of degrees 
$d_1$, $\ldots, d_r$. Assume that $\sum_j d_j \leq n$. Then, if $Z = D_1 
\cap \ldots \cap D_r$, 
\eq{CohInter}{ H^0(Z, W\sO_{Z,K}) = K, \quad\quad H^i(Z, W\sO_{Z,K}) = 0 \ \ 
\mathrm{for}\ i \geq 1. }
\end{prop} 

\begin{proof}
We proceed by induction on $r$. Let $r = 1$, $D = D_1 = Z$, $d = d_1$.
Then the condition $d \leq n$ implies that the homomorphisms
$H^i(\P^n, \sO_{\P^n}) \to H^i(D, \sO_{D})$ are isomorphisms for all
$i$. It follows that, for any $m \geq 1$,
$$ H^0(D, W_m\sO_D) = W_m(k), \quad\quad H^i(D, W_m\sO_D) = 0 \ \ 
\mathrm{for}\ i \geq 1, $$ 
which implies \eqref{CohInter}.

For arbitrary $r$, let $D = D_1 \cup \ldots \cup D_r$, which is a 
hypersurface of degree $d = \sum_j d_j \leq n$. For each 
sequence $1 \leq i_0 < \cdots < i_s \leq r$, let $Z_{i_0,\ldots,i_s} 
= D_{i_0} \cap \cdots \cap D_{i_s}$.  By Corollary \ref{LongMayerViet}, we obtain an
exact sequence 
$$ 0 \!\to\! W\sO_{D,K} \!\to\! \prod_{i=1}^r W\sO_{Z_i,K} \!\to\! 
\cdots \!\to\! \prod_{i=1}^r W\sO_{Z_{1,\ldots,\hat{i},\ldots,r},K} \!\to\! 
W\sO_{Z,K} \!\to\! 0. $$
Applying the previous result to $D$, and the induction hypothesis to
all $Z_{i_0,\ldots,i_j}$ for $j \leq r-2$, we can view this exact
sequence as providing a $\Gamma(\P^n, -)$-acyclic left resolution of
$W\sO_{Z,K}$. Taking sections and observing that the complex
$$ 0 \to K \to \prod_{i=1}^r K \to \cdots \to \prod_{i=1}^r K \to 0 $$ 
is acyclic, except in degree $0$ where its cohomology is equal to 
$K$, we obtain \eqref{CohInter}. 
\end{proof}

\noindent\textit{Proof of Corollary \ref{ThAx}}. Combining 
Proposition \ref{CohAx} with Corollary \ref{zetaW}, we obtain
$$ \zeta^{<1}(Z, t) = \frac{1}{1-t}, $$ 
which, by {\bf } Proposition \ref{Divis}, is equivalent to congruence \eqref{CongAx}
for all finite extensions of $\F_q$. \hfill $\Box$

\medskip
We now discuss some cases where, given a morphism $f : X \to Y$
between two varieties over a finite field, Theorem \ref{SlopeHc}
provides congruences between the numbers of rational points on $X$ and
$Y$.

\begin{prop}\label{IsomWOX}\hspace{-8mm}
\romain Let $f : X \to Y$ be a proper morphism between two separated
$\F_q$-schemes of finite type. If the induced homomorphisms $f^{\ast}
: H^i\cpt(Y, W\sO_{Y,K}) \to H^i\cpt(X, W\sO_{X,K})$ are isomorphisms
for all $i \geq 0$, then 
\eq{CongIsomWOX}{ |X(\F_q)| \equiv |Y(\F_q)| \mod q. } 

\romain Let $X$ be a proper scheme over $\F_q$, and $G$ a finite group
acting on $X$ so that each orbit is contained in an affine open
subset. If the action of $G$ on $H^i(X, W\sO_{X,K})$ is trivial for
all $i$, then
\eq{CongTrivWOX}{ |X(\F_q)| \equiv |(X/G)(\F_q)| \mod q. } 
\end{prop}

\begin{proof}
Assertion (i) is an immediate consequence of \ref{zetaW} and
\ref{Divis}, and implies assertion (ii) thanks to the following lemma:

\begin{lem}\label{WittCohQuot}
Let $X$ be a proper scheme over $\F_q$, $G$ a finite group acting on
$X$ so that each orbit is contained in an affine open subset, $f : X
\to Y := X/G$ the quotient map. Then $f^{\ast}$ induces canonical
isomorphisms
\eq{WCohQuot}{ H^i(Y, W\sO_{Y,K}) \xrightarrow{\ \sim\ \,} 
H^i(X, W\sO_{X,K})^G }
for all $i$.
\end{lem}

The canonical morphism $f : X \to Y$ is finite, and $\sO_Y
\xrightarrow{\ \sim\ \,} f_{\ast}(\sO_X)^G$. By induction on $m$, $R^i
f_{\ast}(W_m\sO_X) = 0$ for all $i \geq 1$ and all $m \geq 1$.
Moreover, the morphisms
$$ W_m\sO_Y \to W_m(f_{\ast}(\sO_X)^G) \to W_m(f_{\ast}(\sO_X))^G 
\cong f_{\ast}(W_m\sO_X)^G$$ 
are isomorphisms for all $m \geq 1$. Since taking invariants under $G$
commutes with inverse limits and with tensorisation by $\Q$, they
provide isomorphisms
$$ W\sO_{Y,K} \xrightarrow{\ \sim\ \,} (f_{\ast}(W\sO_X)^G)_K 
\xrightarrow{\ \sim\ \,} f_{\ast}(W\sO_{X,K})^G. $$
As $\Char(K) = 0$, taking invariants under $G$ commutes with
cohomology for $K[G]$-modules, and we can write 
$$ H^i(Y, W\sO_{Y,K}) \xrightarrow{\ \sim\ \,} 
H^i(Y, f_{\ast}(W\sO_{X,K})^G) \xrightarrow{\ \sim\ \,} 
H^i(Y, f_{\ast}(W\sO_{X,K}))^G. $$ 
On the other hand, the inverse system $(f_{\ast}(W_m\sO_X))_{m \geq
1}$ is $\varprojlim$-acyclic, and $f_{\ast}$ commutes with
tensorisation with $K$. Therefore the morphism
$$ f_{\ast}(W\sO_{X,K}) \to \R f_{\ast}(W\sO_{X,K}) $$ 
is an isomorphism, and we obtain  
$$ H^i(Y, f_{\ast}(W\sO_{X,K}))^G \xrightarrow{\ \sim\ \,} 
H^i(X, W\sO_{X,K})^G, 
$$ 
which proves the lemma.
\end{proof}

\subsection{}\label{ProofIsomOX}
In the most favorable cases, the assumptions of the previous 
proposition can be checked in characteristic $p$. This is the case 
under the assumptions of Corollary \ref{IsomOX}, the proof of which 
follows:
\medskip

\noindent\textit{Proof of Corollary \ref{IsomOX}}.
Let $f : X \to Y$ be a morphism between two proper $\F_q$-schemes such
that the induced homomorphisms $f^{\ast} : H^i(Y, \sO_Y) \to H^i(X,
\sO_X)$ are isomorphisms for all $i \geq 0$. By induction on $m$, it
follows that $H^i(Y, W_m\sO_Y) \to H^i(X, W_m\sO_X)$ is an isomorphism
for all $i \geq 0$ and all $m \geq 1$. As $X$ and $Y$ are proper, 
we obtain that 
$$ H^i(Y, W\sO_{Y,K}) \xrightarrow{\sim} H^i(X, W\sO_{X,K}) $$ 
for all $i \geq 0$. Then congruence \eqref{CongIsomOX} again follows 
from \ref{zetaW} and \ref{Divis}.
\hfill$\Box$

\begin{cor}\label{Primetop}
Let $X$ be a proper scheme over $\F_q$, and $G$ a finite group acting
on $X$ so that each orbit is contained in an affine open subset of
$X$. If $|G|$ is prime to $p$, and if the action of $G$ on $H^i(X,
\sO_X)$ is trivial for all $i$, then
\eq{CongPrimetop}{ |X(\F_q)| \equiv |(X/G)(\F_q)| \mod q. } 
\end{cor}

\begin{proof}
When $|G|$ is prime to $p$, taking invariants under $G$ is an exact 
functor on $k[G]$-modules, and we obtain isomorphisms 
\begin{eqnarray*}
H^i(Y, \sO_Y) \xrightarrow{\ \sim\ \,} 
H^i(Y, (\R f_{\ast}\sO_X)^G) & \xrightarrow{\ \sim\ \,} & 
H^i(Y, (\R f_{\ast}\sO_X))^G \notag \\
& \xrightarrow{\ \sim\ \,} & 
H^i(X, \sO_X)^G = H^i(X, \sO_X). 
\end{eqnarray*} 
Hence we can apply Corollary \ref{IsomOX}.
\end{proof}

\begin{rmks}
Without assumption on $|G|$, Fu and Wan have proved \cite[Theorem
0.1]{WF} that the congruence \eqref{CongTrivWOX} holds under the
following hypotheses:

\alphab $X$ is the reduction of a projective and smooth $W$-scheme 
$X'$ with $p$-torsion free Hodge cohomology;

\alphab The action of $G$ on $X$ is the reduction of an action on 
$X'$ such that the induced action on $H^i(X', \sO_{X'})$ is trivial for 
all $i$. 

We do not know whether, under these hypotheses, the action of $G$ is
trivial on the spaces $H^i(X, W\sO_{X,K})$. However, the following
example seems to indicate that, when $p$ divides $|G|$, a mod $p$
assumption on the action of $G$ as in \ref{Primetop} does not suffice
to provide congruences such as \eqref{CongTrivWOX}.

Let $p = 2$, let $E_1$ be an elliptic curve over $\F_q$ with an
$\F_q$-rational point $t$ of order $2$, and let $E_2$ be another
elliptic curve over $\F_q$. We define $X = E_1 \times E_2$, and we let
the group $G = \Z/2\Z$ act on $X$ via $(x,y) \mapsto (x+t,-y)$, so
that $Y = X/G$ is the classical Igusa surface \cite{Ig}, which is
smooth over $\F_q$. The action of $G$ on $H^i(X, \sO_X)$ is trivial
for all $i$, but, thanks to the K\"unneth formula in crystalline
cohomology, one checks easily that $H^2(X, W\sO_{X,K}) \cong
H^2(X/K)^{<1} \neq 0$, while $H^2(Y,W\sO_{Y,K}) \cong
(H^2(X/K)^G)^{<1} = 0$. Therefore $\zeta^{<1}(Y,t) \neq
\zeta^{<1}(X,t)$, and a congruence such as \eqref{CongTrivWOX} cannot
hold for all powers of $q$.

Note also that, for any $p$, taking the quotient of an abelian variety
by the subgroup generated by a rational point of order $p$ provides an
example where the assumptions of \ref{IsomWOX} (ii) are satisfied,
while those of \ref{IsomOX} fail to be true.
\end{rmks}

\newpage

\bibliographystyle{plain}
\renewcommand\refname{References}

\end{document}